\documentclass{article}
\usepackage{latexsym}
\usepackage{amssymb}
\usepackage{amsmath}
\usepackage{enumerate}
\usepackage[dvips]{graphicx}
\usepackage[latin1]{inputenc}
\title{On cusps  and  flat tops}
\author{Neil Dobbs}

\newcommand\beginpf{\noindent \emph{Proof: }}
\newcommand\eprf{\hfill$\Box$}
\newcommand\epsgamma{cusp }
\newcommand\arr{\mathbb{R}}

\newcommand\dist{\mathrm{dist}}

\newcommand\Crit{\mathrm{Crit}}

\newcommand\scrM{{\cal M}}
\newcommand\scrP{{\cal P}}

\newcommand\scrQ{{\cal Q}}

\newcommand\omu{{\overline{\mu}}}

\newcommand\HD{\mathrm{HD}}

\newcommand\zzz{\mathbb{Z}}
\newcommand\remark{\noindent \emph{Remark: }}
\newcommand\eg{e.g.\ }

\newtheorem{thm}{Theorem}
\newtheorem{dfn}[thm]{Definition}
\newtheorem{lem}[thm]{Lemma}
\newtheorem{prop}[thm]{Proposition}
\newtheorem{cor}[thm]{Corollary}
\newtheorem{conj}[thm]{Conjecture}

\begin{document}
\date{\today}
\maketitle
\begin{abstract}
Non-invertible Pesin theory is developed for a class of piecewise smooth interval maps which may have unbounded derivative, but satisfy a property analogous to $C^{1+\epsilon}$. 
The critical points are not required to verify a non-flatness condition, so the results are applicable to $C^{1+\epsilon}$ maps with flat critical points. If the critical points are too flat, then no 
absolutely continuous invariant probability measure can exist. This generalises a result of Benedicks and Misiurewicz.
\end{abstract}
\section{Introduction}
The ergodic theory of real one-dimensional dynamical systems has been a topic of intense study in recent decades. Especially, much progress has been made for smooth maps of the interval with non-flat critical points (\cite{DeMeloVanStrien, Ledrappier:AbsCnsInterval, BRLSVS, BruinShenVS:NoGrowth} to cite a few), and the theory underlying such dynamical systems is now well-understood. In this work we aim to develop some aspects of the theory beyond the smooth, non-flat setting. The principal results are Theorems \ref{thm:response}, \ref{thm:acimequiv}, \ref{thm:UnstableBis} and \ref{thm:NeilUnstableManif}. 

\subsection{Smooth maps with flat tops}

For smooth maps without a non-flatness condition on the critical points, results are limited. Benedicks and Misiurewicz (\cite{MisBen:Flat}) showed that under a non-recurrence condition, unimodal maps with negative Schwarzian derivative have an 
ergodic, \emph{absolutely continuous, invariant, probability measure} (\emph{acip}) if and only if the logarithm of the derivative is Lebesgue-integrable (see \cite{MeBS08} for the recent exponential family analogue of this result and \cite{Zweimuller:Flat} for statistical properties of the maps considered in \cite{MisBen:Flat}). 
Thunberg, in \cite{Thunberg:Flat}, showed Benedicks-Carleson type results for  unimodal families of maps with critical behaviour like $\exp(-|x|^{-\alpha})$ for $\alpha < 1/8$. He asked whether for $\alpha \geq 1$ no acip can exist.

For maps with \emph{non-flat} critical points (that is, with critical behaviour like $|x|^l$ for some $l > 1$) the log of the derivative is integrable. The maps considered by Thunberg have flat critical points, and the log of the derivative is integrable if and only if $\alpha < 1$. Then by \cite{MisBen:Flat} no acip exists when both $\alpha \geq 1$ and the critical point  is non-recurrent. Non-recurrence is historically a very important condition, but it is quite restrictive, see \cite{Duncan:Rare, Magnus:Rare}. 

For the maps considered in \cite{MisBen:Flat}, Lebesgue measure is ergodic, so any acip necessarily has positive entropy, and thus positive Lyapunov exponent --- see Ruelle's inequality below.

We drop the non-recurrence condition and the negative Schwarzian condition, extending  the \emph{only if} part of Benedicks and Misiurewicz's result to more general  $C^{1+\epsilon}$ maps:
\begin{thm} \label{thm:response}
Let $f: I \to I$ be a $C^{1+\epsilon}$ map of the compact interval $I$ which is piecewise monotone.  Suppose $f$ has an ergodic acip
$\mu$ with positive Lyapunov exponent.
Then the support of $\mu$ is a finite union of intervals on which
\begin{equation}\label{eqn:integrability}
\int \log |Df(x)| dx > -\infty,
\end{equation}
where integration is with respect to Lebesgue measure. 
\end{thm}
\begin{dfn} \label{dfn:mono}
    If $A \subset \arr$, we say a function $g : A \to \arr$ is \emph{piecewise monotone} if and only if there is a finite collection $\scrQ$ of pairwise disjoint intervals, whose union contains $A$, such that for each $Q \in \scrQ$, the restriction $g_{|Q}$ of $g$ to $Q$ is monotone (and not necessarily \emph{strictly} monotone). We call such a $\scrQ$ a \emph{finite partition into intervals of monotonicity}. 
\end{dfn}
Since the derivative is bounded in the smooth setting, $\int \log|Df| d\mu \in [-\infty,+\infty)$, and in particular it exists. \emph{Ruelle's inequality} for $C^1$ maps (\cite{Ruelle:Inequality}) gives
$$h_\mu \leq \max\left(0, \int \log|Df| d\mu\right),$$
where 
  $h_\mu$ denotes the entropy of $\mu$,
so one can replace positive Lyapunov exponent by positive metric entropy in the hypotheses, if one so desires. 
Both our proof and that of \cite{MisBen:Flat} are based on showing that the density of the hypothetical acip is bounded away from zero. In the non-recurrent Misiurewicz setting, one then shows that the first return time to a sufficiently small critical neighbourhood is integrable if and only if (\ref{eqn:integrability}) holds. 
On the other hand, for $C^2$ unimodal maps whose periodic points are all hyperbolic repelling, if the critical point is \emph{recurrent} then the return time to  a critical neighbourhood is always integrable. In particular, the techniques of \cite{MisBen:Flat} do not work in the recurrent general setting.

The following corollary  responds to  the question of Thunberg.
\begin{cor}\label{cor:acimaone}
 Let $g_b$ be a map from the unimodal family 
$$
x \mapsto -1 + b\left(1 - e^{-1 - |x|^{-\alpha}} \right).
$$
If $\alpha \geq 1$ then no ergodic acip with positive Lyapunov exponent exists.
\end{cor}
\beginpf
Were such a measure to exist,  its support would contain an interval, by Theorem \ref{thm:response}, and so would necessarily contain the critical point at 0.
Then Theorem \ref{thm:response} would imply that $\log|Dg_b(x)|$ is integrable with respect to Lebesgue measure.
But $\log|Dg_b(x)| = h(x) - |x|^{-\alpha}$, where $h$ is some function integrable with respect to Lebesgue measure. If $\alpha \geq 1$ then $|x|^{-\alpha}$ is not integrable.
\eprf

\subsubsection{On hypotheses and limitations}
For these results, ergodicity is not strictly necessary. As per \cite{Ledrappier:AbsCnsInterval}, if one just supposes that $\mu$ is absolutely continuous and that $\mu$ almost every point has positive (pointwise) Lyapunov exponent, then one can use the ergodic decomposition theorem and continue similarly. In the interest of brevity and clarity we do not do this here. 

The assumption that the acip has positive Lyapunov exponent is not wholly unnatural.
If there exists a $\delta >0$ such that almost every point is contained in arbitrarily small intervals each mapped with universally bounded distortion by some iterate of $f$ onto some interval of length $\delta$, then it is not hard to show (with a density point argument) that any absolutely continuous measure actually has positive metric entropy, and that there are at most a finite number of such measures. Then, by Ruelle's inequality, the measure has positive Lyapunov exponent. Let us say that such maps have the \emph{large scale property}.

In the non-flat setting, it is known that Lebesgue measure is ergodic \cite{Martens:Cantor, BlokhLyubich}, and that any acip must have positive entropy. Outside the non-flat setting (and in particular for the unimodal family considered here), this is unknown and difficult. There may conceivably be maps with acips with negative or zero Lyapunov exponents. However, it should be possible to construct maps from the above unimodal family which have the large scale property, even if $\alpha \geq 1$, by controlling the rate of recurrence of the critical point. For such maps, Lebesgue measure would be conservative and ergodic, but no acip would exist. 

While such maps probably exist, we do not believe that there is a positive measure set of parameters for which Lebesgue measure is conservative. We call a unimodal map \emph{hyperbolic} if it has an attracting periodic orbit with a  basin of attraction whose complement is (uniformly) hyperbolic repelling. For hyperbolic unimodal maps, Lebesgue almost every point lies in the basin of attraction, and Lebesgue measure is not conservative.
\begin{conj} Consider the unimodal family 
    $$g_b : x \mapsto -1 + b\left(1 - e^{-1 - |x|^{-\alpha}} \right).   $$
If $\alpha \geq 1$, then for Lebesgue almost every $b$, the map $g_b$ is hyperbolic.
\end{conj}

In the theorem, we suppose that $f$ is piecewise monotone. We need this hypothesis to find a finite generating partition (see Proposition \ref{prop:partition}). Since the partition we find is finite, its entropy is finite. Without assuming finiteness of the number of turning points, it may be possible to find an infinite generating partition with finite entropy and other good properties, but we have not succeeded in doing so. Note that we do not exclude $f$ having an infinity of inflection points, or, \emph{a priori}, $f$ having an interval of critical points, so it is a little more general than assuming that the set of critical points is finite.

\subsection{Maps with unbounded derivative}

The second goal of the paper is to develop non-invertible Pesin theory for a class of  maps with discontinuities and unbounded derivative. 
Map with unbounded derivative are of interest due to their links with the Lorenz map. 
See \cite{LuzzattoViana:Ast00} and \cite{LuzzattoTucker} for a discussion of this and \cite{Luzzatto:Lorenz} and \cite{DHL06} for existence results for absolutely continuous, invariant, probability measures. An early, general and very powerful result is by Rychlik \cite{Rychlik:BV}. 

In the following section  we shall introduce a new class of maps called cusp maps. This class of maps will include all $C^{1+\epsilon}$ maps. 
For piecewise $C^{1+\epsilon}$ maps whose critical points verify a non-flatness condition, Pesin theory was studied by Ledrappier in \cite{Ledrappier:AbsCnsInterval}. Given a measure with positive Lyapunov exponent, he showed existence of the  unstable manifold in the natural extension, and several results which follow from it. The non-flatness hypothesis, used repeatedly in his proofs,  means that we cannot use his results in the proof of Theorem \ref{thm:response}.

In Theorem \ref{thm:NeilUnstableManif}, whose statement is overly technical for this introduction, we show existence of the unstable manifold for cusp maps. Even for $C^{1+\epsilon}$ maps our result is stronger than that of Ledrappier since we do not assume non-flatness of critical points. Moreover our proof is more direct. For $C^{1+\epsilon}$ maps one can also, with some work, deduce this result from \cite{Newhouse:EntropyVolume}; however the proof in that higher-dimensional setting is considerably more complex. 
%We are not aware of anyone using \cite{Newhouse:EntropyVolume} in one-dimensional dynamics. 
This will allow us to prove Theorem \ref{thm:response}.

In \cite{Me:MaxEntLMS}, the author gave a $C^r$ version of Ledrappier's unstable manifold theorem and used it to prove $C^r$ conjugacy results. We shall also state a $C^r$ version here, but shall refer to \cite{Me:MaxEntLMS} for the proof.

With unstable manifold in hand, we use regularly returning (or nice) intervals to give simple proofs of the dynamical volume lemma in Proposition \ref{prop:dvl} and of the existence of a Pesin partition in Proposition \ref{thm:pesin}, compare \cite{HofbauerRaith:HDmeasure} and \cite{Ledrappier:AbsCnsInterval}. 

Given a transformation $g$, we denote by $\scrM(g)$ the collection of ergodic $g$-invariant probability measures. If $f$ is a cusp map, $\mu \in \scrM(f)$ and $\chi_\mu := \int \log |Df| d\mu$ exists,  then we call $\chi_\mu$ the Lyapunov exponent of $\mu$. For cusp maps, $\chi_\mu$ could be $+\infty$ or (which we do \emph{not} exclude) $-\infty$, or it need not exist at all. In Section \ref{sec:cheb} there is an example of a cusp map with a measure of maximal entropy with Lyapunov exponent $+\infty$.

The paper culminates with the proof of the following result.

\begin{thm}\label{thm:acimequiv}
    Let $f$ be a \epsgamma  map which is piecewise monotone. Suppose $\mu$ is an ergodic, invariant, probability measure for $f$ with positive finite Lyapunov exponent $\chi_\mu$.
The following conditions are equivalent:
\begin{enumerate}
\item \label{enum:ace1}$\mu$ is absolutely continuous with respect to Lebesgue measure;
\item \label{enum:ace2}$h_\mu = \chi_\mu$;
\item \label{enum:ace3}the density of $\mu$ with respect to Lebesgue measure  is bounded from below by a positive constant on an open interval;
\item \label{enum:ace4}$\mu$ is generated by a full expanding induced Markov map  with integrable return time.
\end{enumerate}
\end{thm}
\beginpf Follows immediately from Corollary \ref{cor:cuspacim} and  Propositions \ref{thm:cuspacim} and \ref{thm:cuspacimMarkov}. 
\eprf

We refer to Section~\ref{sec:markov} for the definition of a full expanding induced Markov map. The set on which the density is bounded away from zero is often large, see Lemma \ref{lem:InvariantDensity} and Theorem \ref{thm:densitybddbelow}.

Ledrappier \cite{Ledrappier:AbsCnsInterval} showed equivalence between \ref{enum:ace1} and \ref{enum:ace2} for $C^{1+\epsilon}$ maps with non-flat critical points. Bruin \cite{Bruin:Markov0} showed equivalence of all four conditions in the case of unimodal maps with non-flat critical points and negative Schwarzian derivative. There is a recent related result in \cite{BT:EquilibriumInterval} for multimodal maps with non-flat critical points and negative Schwarzian derivative. Theorem \ref{thm:acimequiv} represents a substantial improvement. 

We assume existence and finiteness of the Lyapunov exponent. As we have already stated, it need not exist. In Proposition \ref{prop:examples}, we construct benign examples of cusp maps with acips for which the Lyapunov exponent does not exist. We also show that cusp maps can have acips despite the presence of smooth parabolic fixed points, and that there are cusp maps for which the measure of maximal entropy has positive but infinite Lyapunov exponent. 

 The only reference measure we consider is Lebesgue measure. In the setting of holomorphic dynamics in \cite{Me:Rations}, we consider  more general conformal measures and invariant measures absolutely continuous with respect to them. It would be interesting to have similar results for conformal measures in the interval setting; we do not attempt this here. 

\subsubsection{An alternative hypothesis: Positive entropy}
There is an alternative to assuming that the derivative should tend to $0$ or $\pm \infty$ at the boundary of the domain of definition of the cusp map. We use the distortion bound and this assumption to guarantee that, when we pull back a small enough ball along a typical branch, we do not meet a discontinuity. An alternative approach is to assume positive entropy (as well as positive, finite, Lyapunov exponent). 

Assuming positive entropy, we  extend our results to a broader class of maps with a slightly more restrictive hypothesis on the measure.
The Theorem \ref{thm:NeilUnstableManif} referred to in the following result says that along most inverse branches there is an interval that can be pulled back with good distortion control

\begin{thm} \label{thm:UnstableBis}
    Let $f: \bigcup_{j=1}^d I_j \to I$ be a piecewise monotone map, a $C^1$ diffeomorphism on each subinterval $I_j$ of the compact interval $I$, which satisfies the conditions (\ref{eqn:epsgamma}) and (\ref{eqn:epsgamma2}). Suppose $\mu \in \scrM(f)$ has positive, finite Lyapunov exponent and positive entropy. 

Then the conclusions of 
Theorem 
\ref{thm:NeilUnstableManif} hold.
\end{thm}

Since it was only in the proof of Theorem \ref{thm:NeilUnstableManif} that the assumption on the derivative at the boundary of the domain of definition was used, the other results also hold:
\begin{thm} 
    Let $f: \bigcup_{j=1}^d I_j \to I$ be a piecewise monotone map, a $C^1$ diffeomorphism on each $I_j$, satisfying the conditions (\ref{eqn:epsgamma}) and (\ref{eqn:epsgamma2}). Suppose $\mu \in \scrM(f)$ has positive, finite Lyapunov exponent and positive entropy. 

Then the conclusions of Theorem \ref{thm:acimequiv} 
and Propositions \ref{prop:partition}, \ref{prop:dvl} and \ref{thm:pesin} hold. 
\end{thm}

The assumption of positive entropy allows one to use the lift of the measure in the Hofbauer tower (or Markov extension) to study recurrence in the natural extension, and in particular to prove a result, Proposition \ref{prop:fibres}, about existence of unstable fibres.  

Proposition \ref{prop:fibres} is a topological (and metric) result. It does not give any information on distortion or lack thereof. For that we use positivity of the Lyapunov exponent. 

For piecewise-monotone maps with bounded derivative, preceding results include \cite{Hofbauer:IntrinsicErgodicityI, Hofbauer:IntrinsicErgodicityII, HofbauerRaith:DynamicVolumeLemma, Keller:density, Keller:Liftability}.

 \subsection{Structure of the paper}
The paper is structured as follows. In the next section we introduce the class of cusp maps. Then we prove some elementary distortion estimates. In Section~\ref{sec:unst} we define the natural extension and prove existence of the unstable manifold. 
In Section~\ref{sec:idvl} we show the existence of useful regularly returning intervals. We then use these intervals to find a finite generating partition with good Markov properties and prove the dynamical volume lemma in Section~\ref{sec:dvl}. In Section~\ref{sec:pesin} we show the existence of a Pesin partition of the natural extension. All results up to here are for an arbitrary ergodic invariant probability measure with positive finite Lyapunov exponent. In Section~\ref{sec:ledrappartition} we start the study of absolutely continuous measures and prove Theorem~\ref{thm:response}. In the following section we study induced Markov maps. Then in Section~\ref{sec:posentropy} we outline the proofs of Theorem~\ref{thm:UnstableBis} and Proposition~\ref{prop:fibres}. To finish up, in Section~\ref{sec:cheb} we present some cusp maps conjugate to the Chebyshev map $x \mapsto 4x(1-x)$ and look at some of their invariant measures.

\begin{section}{Definition of cusp maps}
%Smooth maps of the interval have generated a lot of interest and a huge body of work \cite{DeMeloVanStrien}. Maps which are not smooth due to the presence of singularities, or poles, are also worth studying due, for example,  to their links with the Lorenz map \cite{LuzzattoTucker}. 
%However, the presence of poles has prevented much progress from being made. At critical points one has negative Schwarzian derivative which allows one to control the accumulation of distortion of iterates. 
%At poles, conversely, the Schwarzian derivative is unboundedly positive. With poles, the fact that an iterate maps an interval diffeomorphically onto nested intervals does not imply that the pullback of the inner interval is relatively far from singular  points.

%See \cite{LuzzattoTucker} and the references therewithin for more information. 
Our goal is to develop the ergodic theory for piecewise smooth interval maps 
 with singularities where the derivative, on at least one side, may tend to infinity. %We shall consider ergodic invariant probability measures with positive \emph{finite} Lyapunov exponent and 

For continuous maps with two smooth monotone branches, if the norm of the derivative tends to infinity as one approaches the turning point, the turning point is called a \emph{cusp}. This leads us to introduce the following definitions. Throughout the paper, we shall denote by $I$ a non-degenerate compact interval.
\begin{dfn}
A map $f: \bigcup_{j} I_j \to I$ is a \emph{basic cusp map} if 
$(I_j)_{j}$ is  a finite or countable 
collection of pairwise disjoint open subintervals of $I$ such that $f$ is a $C^1$ diffeomorphism from each interval $(p_j,q_j) := I_j$ onto its image $f(I_j)$
and such that 
the following limits exist and equal either 0 or $\pm \infty$:
$$
\lim_{x \to p_j^+} Df(x), \lim_{x \to q_j^-} Df(x).
$$
\end{dfn}
Notice that $f$ is defined on a union of \emph{open} intervals, a fact that is used repeatedly throughout the paper. For example, it implies that, given an interval $U$, $f^n$ is a diffeomorphism on each connected component of $f^{-n}(U)$, and each of these connected components is open. Choosing intervals to be open in the definition simplifies life later on, and from a measure-theoretical point of view, not having boundary points is unimportant. 

We also note that the domain of definition of $f$ may have gaps, that is, the closure of the domain of definition need not be the entire interval $I$.

%The set $X := I \setminus \bigcup_j I_j$ is called the \emph{singular set} or the \emph{set of singularities of $f$}. 
%%
%Let the \emph{invariant set $K_f$} be the set of all $x$ such that $f^n(x) \in I\setminus X$ for all $n \geq 0$. It is forward-invariant but generally is not compact.
Let $A \subset \arr$. Recall that a differentiable map $g:A \to \arr$ is $C^{1+\epsilon}$ with constants $C$ and $\epsilon$ if, for all $x,x'$ in $A$, 
$$
|Dg(x) - Dg(x')| \leq C |x-x'|^\epsilon.
$$
%If $f : \bigcup I_j \to I$ is a basic cusp map, then each restriction $f_{|I_j} : I_j \to I$ extends continuously to a map $f_j : \overline{I_j} \to I$ and, if one compactifies $\arr$ to get $\overline{\arr} = \arr \cup \{+ \infty\} \cup \{-\infty\}$, then $Df_j : \overline{I_j} \to \overline{\arr}$ is continuous.
\begin{dfn} \label{dfn:cursep}
    A basic cusp map $f: \bigcup_j I_j \to I$ is a \emph{\epsgamma map} (\emph{with constants $C, \epsilon$}) if there exist constants  $C, \epsilon > 0$ such that on each $I_j$, 
    %$\overline{I_j}$,
    \begin{itemize}
    \item for all $x,x'$ such that $|Df_j(x)|, |Df_j(x')| \leq 2$,
    \begin{equation}\label{eqn:epsgamma}
        |Df_j(x) - Df_j(x')| \leq C |x-x'|^\epsilon;
    \end{equation}
    \item for all $x,x'$ such that $|Df_j(x)|, |Df_j(x')| \geq 2^{-1}$,
    \begin{equation}\label{eqn:epsgamma2}
        \left|\frac{1}{Df_j(x)} - \frac{1}{Df_j(x')}\right| \leq C |x-x'|^\epsilon.
    \end{equation}
    \end{itemize}
\end{dfn}
One can check that poles of the form $x^\alpha$ where $0 < \alpha < 1$ (these are \emph{poles of root type}) satisfy this latter relation. In figure \ref{fig:somecuspmaps} we present (the graphs of) some cusp maps.
\begin{figure}[htb]
\center{\includegraphics[totalheight=0.4\textheight]{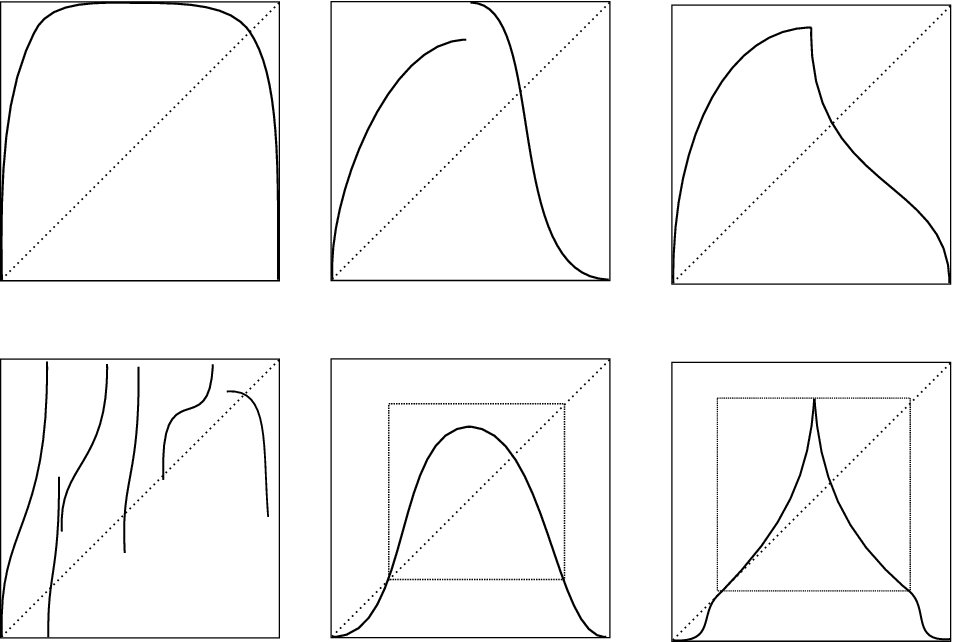}}
%\center{\includegraphics{../../Images/somecuspmaps}}
\caption{Some cusp maps. \label{fig:somecuspmaps}}
\end{figure}

\begin{dfn} A continuous map $g : I \to I$ which is $C^1$ on all but a countable set $S$ and for which $g_{|I\setminus(S\cup \Crit)}$ is a cusp map will also be called a cusp map. We use $\Crit$ to denote the set of critical points of $g_{|I\setminus S}$. 
\end{dfn}
Since all our results  will be for ergodic probability measures with positive finite Lyapunov exponent, the distinction between the continuous map and its restriction to $I\setminus(S\cup \Crit)$ will be of no importance.

We are also interested in some analogue to $C^r$ maps. We want a condition which is satisfied by higher derivatives $D^if$ of $f$ and which holds for all poles of root type. What appeared naturally and suffices to prove nice $C^r$ distortion properties for induced Markov maps is the following. 
\begin{dfn} Let $r \in \{2, 3, \ldots\}$. A \epsgamma map $f : \bigcup_j I_j \to I$ is a \emph{$C^r$ cusp map} (\emph{with constants $p,C$}) if there exist constants $C, p > 1$ such that, for each $j$ one has: 
\begin{itemize}
\item
$f_j$ is $C^r$ on $I_j$;
\item for all $x \in I_j$ such that $0 < |Df_j(x)| \leq 2$, and for all $i$ such that $2 \leq i \leq r$, $|D^if(x)| < C$;
\item for all $x \in I_j$ with $|Df_j(x)| \geq 2$, and for all $i$ such that $2 \leq i \leq r$, 
\begin{equation}\label{eqn:crcusp}
\frac{|D^if(x)|}{|Df(x)|^p} < C;
\end{equation}
\end{itemize}
If $f$ is a $C^r$ cusp map for all integers $r > 1$ then $f$ is a $C^\infty$ cusp map.
\end{dfn}
%\remark If a map $f : I \to I$ coincides with a restriction of a cusp map on $I \setminus \Crit(f)$ then we may also call $f$ a cusp map. For example, a map $x \mapsto ax(1-x)$ on $[0,1]$ is a cusp map. 
\begin{dfn}
Let $f : I \to I$. A \emph{pole (at $p$)} is an interval $(p, p+\epsilon) \subset I$ such that $|\varepsilon| > 0$,  $f$ is $C^1$ on $(p, p+\epsilon)$ and
$$
\lim_{h \to 0^+} |Df(p + h\epsilon)| = \infty.
$$
\end{dfn}
%\begin{dfn} \label{dfn:indeptsing}
%A cusp-map $f$ has \emph{independent singularities} if the singular set $X$ satisfies $f^n(\partial X) \cap \partial X = \emptyset$ for all $n > 0$.
%\end{dfn}

% for commonly studied measures (measures of maximal entropy, absolutely continuous measures). Finiteness of the integral is crucial to control approach rates to singularities. 

\end{section}
\begin{section}{Distortion estimates}\label{sec:distnest}
We deduce some simple distortion estimates, culminating in Lemma \ref{lem:distnkey} which roughly speaking says that the distortion is well-bounded on a ball of radius depending polynomially on the derivative.

In the following two lemmas, one can think of $\phi$ as $Df$ if $Df$ is small, or as $\frac{1}{Df}$ if $Df$ is big.
\begin{lem}\label{lem:oneepsdistn}
    Let $c, C, \epsilon >0$ be  positive constants.
Let $\phi: A \to \arr$ satisfy $|\phi(x) - \phi(x')| \leq C|x-x'|^\epsilon$ for all $x, x' \in A$.

Then, if $x,x' \in A$  and $|\phi(x)| > c$,
\begin{equation}\label{eqn:oneepsdistn}
1 - Cc^2 \frac{|x-x'|^\epsilon}{c^3} \leq \frac{\phi(x')}{\phi(x)} \leq 1+Cc^2\frac{|x-x'|^\epsilon}{c^3}.
\end{equation}
\end{lem}
\beginpf
We have $|\phi(x) - \phi(x')| \leq C|x-x'|^\epsilon$ and $|\phi(x)| > c$. Thus
$$
\frac{|\phi(x) - \phi(x')|}{|\phi(x)|} \leq  C \frac{|x-x'|^\epsilon}{c} = C c^2\frac{|x-x'|^\epsilon}{c^3}.
$$
One recovers (\ref{eqn:oneepsdistn}) upon rewriting the left-hand side as $\left|1 - \frac{\phi(x')}{\phi(x)}\right|$ and using the triangle inequalities.
\eprf
\begin{lem}\label{lem:oneepsdistnlog}
Let $C,\epsilon > 0$ be constants.
There exists a constant $c_0$, depending on $C, \epsilon$ and satisfying $0 <c_0 < \log 2$, with the following property:

If $\phi: A \to \arr$ satisfies $|\phi(x) - \phi(x')| \leq C|x-x'|^\epsilon$ for all $x, x' \in A$,
and if $c$  satisfies $0 < c \leq c_0$, 
then for all $x,x' \in A$ such that $|x-x'|^\epsilon < c^3$ and $|\phi(x)| > c$,
\begin{equation}\label{eqn:oneepsdistnlog}
\left|\log |\phi(x)| - \log |\phi(x')| \right| \leq c \frac{|x-x'|^\epsilon}{c^3} < c.
\end{equation}
\end{lem}
\beginpf
There exists a constant $a_0$, with $0 < a_0 < \log 2$, such that $\log \left(1-\frac{a}{2}\right) > -a$ for all $a$ satisfying $0 < a < a_0$. Choose $c_0$ with $0 < c_0 < a_0$ and  $Cc_0 < \frac{1}{2}$. Then by Lemma \ref{lem:oneepsdistn}, for $c,x$ such that $0<c\leq c_0$ and $|\phi(x)|>c$, 
$$ 
1 - \frac{c}{2}\frac{|x-x'|^\epsilon}{c^3} \leq 1 - Cc^2\frac{|x-x'|^\epsilon}{c^3} \leq \frac{\phi(x')}{\phi(x)} \leq 1+Cc^2\frac{|x-x'|^\epsilon}{c^3} \leq 1 + \frac{c}{2}\frac{|x-x'|^\epsilon}{c^3}.
$$
Taking logs and using $c \leq c_0 < a_0$,
$$
-c\frac{|x-x'|^\epsilon}{c^3} \leq \log\left(\frac{\phi(x')}{\phi(x)}\right) \leq c\frac{|x-x'|^\epsilon}{c^3},
$$
as required.
\eprf
\begin{lem}\label{lem:distnkey}
Let $f : \bigcup_j I_j \to I$ be a \epsgamma map with constants $C,\epsilon$. 

There exists a constant $c_0 >0$ such that if $c$  satisfies $0 < c < c_0$ and 
 $x,x' \in \bigcup_j I_j$ satisfy $|x-x'|^\epsilon < c^3$ and $c < |Df(x)| < c^{-1}$, then $x$ and $x'$ are in the same component of the domain of definition of $f$ and
\begin{equation}\label{eqn:oeedistnlog}
\left|\log |Df(x)| - \log |Df(x')| \right| \leq c \frac{|x-x'|^\epsilon}{c^3} < c.
\end{equation}
\end{lem}
\beginpf
Let $c_0$ be given by Lemma \ref{lem:oneepsdistnlog}.  
Fix $x \in \bigcup_j I_j$ and $c$ with $0 < c< c_0$, such that $c < |Df(x)| < c^{-1}$. Let $A$ denote the connected component of 
$$\{x' \in B(x,c^{(3/\epsilon)}) : 2^{-1}|Df(x)| < |Df(x')| < 2|Df(x)|\}
$$ 
containing $x$; in particular $A$ is contained in a connected component of the domain of definition of $f$. 
If $|Df(x)| \leq 1$, set $\phi := |Df|$, otherwise set $\phi := |Df|^{-1}$. By Lemma \ref{lem:oneepsdistnlog}, we have   
$$|\log |Df(x)| - \log |Df(x')|| \leq c \frac{|x-x'|^\epsilon}{c^3} <c < \log 2$$
for all $x' \in A$. It follows immediately that  $A =B(x,c^{(3/\epsilon)})$.
\eprf
\end{section}

%\begin{section}{The natural extension}\label{sec:ext}
\begin{section}{Unstable Manifold} \label{sec:unst}
%%%%%%% REDEFINE THIS AS PER THE RATIONAL SECTION? %%%%%%%%%%%%
We define the natural extension as per \cite{Ledrappier:AbsCnsComplex}. Let 
$$
Y := \{ y = (y_0y_1y_2\ldots) : f(y_{i+1}) = y_i \in I\}.
$$
Define $F^{-1} : Y \to Y$ by $F^{-1}((y_0y_1\ldots)) := (y_1y_2\ldots)$. Then $F^{-1}$ is invertible with inverse $F : F^{-1}(Y) \to Y$. The projection $\Pi : Y \to I$ is defined by $\Pi :y = (y_0y_1\ldots)  \mapsto y_0$. 
Then $\Pi \circ F = f \circ \Pi$. Given any measure $\mu \in \scrM(f)$ there exists a unique $F$-invariant measure $\omu$ such that $\Pi_*\omu = \mu$. Moreover $\omu \in \scrM(F)$ and $\omu \in \scrM(F^{-1})$ (see \cite{Rohlin:Exact}). 

We call the triplet $(Y,F,\omu)$ the \emph{natural extension} of $(f,\mu)$ (it is also called the Rohlin extension or the canonical extension). 

Let us remark that invariant probability measures give no mass to the sets of points $x$ for which there is an $n > 0$ such that $f^n(x)$ is not defined, nor do they give mass to the set of $x$ for which there exists an $n>0$ and \emph{no solution} $x'$ to $f^n(x') = x$. Thus, $F^n(y)$ is defined for all $n \in \zzz$ for $\omu$ almost every $y \in Y$.
\begin{thm}\label{thm:NeilUnstableManif}
Let $f : \bigcup_j I_j \to I$ be a \epsgamma map. Suppose $\mu \in \scrM(f)$ has positive finite Lyapunov exponent $\chi = \int \log|Df| \, d\mu > 0.$ Denote by $(Y, F, \overline{\mu})$ the natural extension of $(f, \mu)$. 

Then there exists a measurable function $\alpha$ on $Y$, $0 < \alpha < 1/2$ almost everywhere, such that for $\overline{\mu}$ almost every $y \in Y$ 
there exists a  set $V_y \subset Y$ with the following properties:
\begin{itemize}
\item $y \in V_y$ and $\Pi V_y =  B(\Pi y, \alpha(y))$;
\item for each $n >0$, $f^n : \Pi F^{-n}V_y \to \Pi V_y$ is a diffeomorphism (in particular it is onto);
\item for all $y' \in V_y$
$$
  \sum_{i=1}^\infty \left| \log|Df(\Pi F^{-i}y')| - \log |Df(\Pi F^{-i}y)|\right|  
< \log 2;
$$
%\begin{equation}\label{eqn:log4b}
%< \log 4 \frac{|\Pi y - \Pi y'|}{|\Pi V_y|} < \log 4;
%\end{equation}
\item for each $\eta >0$ there exists a measurable function $\rho$ on $Y$, $0< \rho(y) < \infty$ almost everywhere, such that 
$$
\rho(y)^{-1}e^{n(\chi - \eta)} < |Df^n(\Pi F^{-n}y)| < \rho(y)e^{n(\chi + \eta)}.
$$
\end{itemize}
In particular, $|\Pi F^{-n}V_y| \leq \rho(y)e^{-n(\chi-\eta)}$.
\end{thm}
The above theorem is sufficient for the purposes of this paper. A stronger version is possible: 
\begin{thm} \label{thm:Unstable2}
Supplementarily to Theorem \ref{thm:NeilUnstableManif}, there exists measurable $\gamma_1 < \infty$ $\overline{\mu}$-almost everywhere such that for all $y' \in V_y$,
$$
 \sum_{i=1}^\infty \left| \log|Df(\Pi F^{-i}y')| - \log |Df(\Pi F^{-i}y)|\right|  
\leq \gamma_1(y) |\Pi y - \Pi y'|^\epsilon.
$$
If $f$ is also a $C^r$ cusp map for some integer $r \geq 2$, then there exists a measurable function $\gamma_r,$   with $\gamma_r(y) < \infty$ $\overline{\mu}$-almost everywhere, such that for all $n$ and all $i = 1,2,\ldots, r-1$,
$$
\frac{1}{|Df^n(\Pi f^{-n}y)|} \left| D^i \log |Df^n(\Pi f^{-n}y)| \right| \leq \gamma_r(y).
$$
\end{thm} 
\beginpf We refer to Proposition 2.11 of \cite{Me:MaxEntLMS} for the proof. Our definition of $C^r$ cusp maps allows that proof to be applied directly. The full proof for cusp maps is also available in \cite{Me:PHD}.
\eprf

The rest of this section is devoted to the proof of Theorem \ref{thm:NeilUnstableManif}. It will be broken up into several lemmas. The strategy is as follows. We have shown in Lemma \ref{lem:distnkey} that we have a good distortion bound on $$ B\left(x,\min(c_0^{3/\epsilon},|Df(x)|^{\pm 3/\epsilon})\right)
$$ 
for some constant $c_0$. Next we show that the derivative $|Df(\Pi F^{-n}y)|$ along backwards orbits is bounded from below by a sub-exponential sequence almost everywhere.  This allows us to define a slowly-shrinking sequence of balls on which one has (slow-) exponentially good distortion. Positive Lyapunov exponent will then imply that the pullbacks of some small ball will always land inside the sequence of balls with exponentially good distortion bounds, so the total distortion will be summable.
 
\noindent \beginpf We will need to swallow up some constants.
Fix $\delta > 0$ such that 
$\delta < \eta$,  
$(\chi - 3\delta) > 3\delta /\epsilon$.
Subsequently fix $N > 0$ large enough that, for all $n \geq N$, the following inequalities hold
\begin{eqnarray}
\label{eqn:sumndelta} 2^{-1} \log 2 +\sum_{m\geq N}  e^{-m\delta} &<& \log 2;\\
e^{-n\delta} &<& c_0;\\
\label{eqn:balls} 
2 e^{-n(\chi-\delta)} &<& 2^{-1} e^{-(n+1)\delta(1 + 3/\epsilon)},
\end{eqnarray}
where $c_0$ comes from Lemma \ref{lem:distnkey}.

\begin{lem} \label{lem:sea} For $\omu$ almost every $y$, there exists $n(y) \geq N$  such that for all $n \geq n(y)$,
$$
2 e^{-n\delta} < |Df(\Pi F^{-n}y)| <  e^{n\delta}/2
$$
and 
$$
 e^{n(\chi - \delta)} \leq |Df^n(\Pi F^{-n}y)| \leq e^{n(\chi + \delta)}.$$
\end{lem}
\beginpf The first holds because the limit of $(1/n)\log |Df^n(\Pi F^{-n}y)|$ exists for almost all $y$; the second because it equals $\chi$. 
\eprf

\begin{lem} \label{lem:bndistn}
 Let $B_n := B(\Pi F^{-n}y, 2^{-1} e^{-n3\delta/\epsilon})$. For all $n \geq n(y)$, 
$B_n$ is contained in the domain of definition of $f$
and, for all $x,x' \in B_n$,
$$
\left| \log|Df(x)| - \log|Df(x')|\right| < e^{-n\delta}.$$
% \frac{\sigma(x,x')}{e^{-n3\delta}}.$$
\end{lem}
\beginpf
Follows from Lemmas \ref{lem:distnkey}, \ref{lem:sea}. 
\eprf
\begin{lem}  For $n \geq n(y)$, $f(B_{n+1}) \supset B(\Pi F^{-n}y, 2 e^{-n(\chi - \delta)})$.
\end{lem}
\beginpf By the preceding lemmas, $|Df(x)| > e^{-(n+1)\delta}$ on $B_{n+1}$, so $f(B_{n+1}) \supset B(\Pi F^{-n}y, 2^{-1}e^{-(n+1)3\delta/\epsilon}e^{-(n+1)\delta})$. Then use (\ref{eqn:balls}).
\eprf

\begin{lem} \label{lem:ballinduct}
 Suppose $n \geq n(y)$ and $V$ is an open ball containing $\Pi y$ with $|V| < 1$ and suppose $V_n$ is such that $V_n \ni \Pi F^{-n}y$ and $f^n: V_n \to V$ is a diffeomorphism with distortion bounded by some $r$ with $0 < r < \log 2$, i.e.,
$$
\left| \log |Df^n(x)| - \log |Df^n(x')|\right| < r < \log 2$$
for all $x,x' \in V_n$. 
Then there exists $V_{n+1} \ni \Pi F^{-(n+1)}y$ such that the map $f^{n+1} : V_{n+1} \to V$ is a diffeomorphism with distortion bounded by 
$$
r + e^{-(n+1)\delta}.$$ 
% \frac{\sigma(x,x')}{e^{-(n+1)3\delta}}$$.
\end{lem}
\beginpf 
We have that $|Df^n(\Pi F^{-n}y)| > e^{n(\chi - \delta)}$ so, by the distortion bound on $f^n$,  $|V_n| < 2e^{-n(\chi - \delta)}$ and $V_n \subset f(B_{n+1})$. The result follows.
\eprf

Now let $V$ be a sufficiently small ball centred on $\Pi y$ that there exists a set $V_{n(y)} \ni \Pi F^{-n(y)}(y)$ such that $f^{n(y)} : V_{n(y)} \to V$ is a diffeomorphism and for all $x,x' \in V_{n(y)}$, 
$$
\sum_{i=0}^{n(y)-1} \left| \log |Df(f^i(x))| - \log |Df(f^i(x'))|\right| \leq (1/2) \log 2.$$
% \frac{\sigma(f^{n(y)}(x), f^{n(y)}(x'))}{|V|}.$$
For $0 \leq n < n(y)$ define $V_n := f^{n(y)-n}(V_{n(y)})$. For $n > n(y)$, define $V_n$ inductively using Lemma \ref{lem:ballinduct}.
For any $n > 0$, for any $x,x' \in V_n$, we have
$$
\sum_{i=0}^{n-1} \left| \log |Df(f^i(x))| - \log |Df(f^i(x'))|\right| \leq (1/2) \log 2 + \sum_{j =n(y)}^\infty e^{-n\delta} < \log 2.$$
Define $V_y$ as the set of $y' \in Y$ such that $\Pi F^{-n}(y') \in V_n$ for all $n \geq 0$. Let $\alpha(y)$ be the radius of $V$. 
    Note that one can choose [$n(y)$ minimal, so it is measurable, and then] $V$ [maximal] so that $\alpha$ is measurable. The existence of $\rho$ is easy, by ergodicity. This completes the proof of Theorem \ref{thm:NeilUnstableManif}.
\eprf

\end{section}

\begin{section}{Regularly returning intervals}  \label{sec:idvl}
The following lemma is simple and known. We include the proof for completeness and as an introduction to the arguments we will use later on, concerning points `going to the large scale'.
\begin{lem} \label{thm:periodicpointsdense}
Let $f$ be a \epsgamma map and suppose $\mu \in \scrM(f)$ has positive finite Lyapunov exponent $\chi_\mu$.

Inside any interval $U$ of positive measure, there are positive measure intervals $W \subset V$ and $n$ large such that $f^n$ maps $W$ diffeomorphically onto $V$ with $|Df^n_{|W}| >2$. 
The closure of the set of repelling periodic points has full ($\mu$-) measure. 
The closure of the collection of inverse images (i.e. the closure of the backward orbit) of some periodic point has full measure.
\end{lem}
\beginpf
Let $U$ be an open interval of positive measure. %For the first part it suffices to show that there are periodic points in  $U$.
Let $(Y, F, \overline{\mu})$ be the natural extension for $(f, \mu)$. Modulo a set of $\overline{\mu}$-measure zero one can write
$$
\Pi^{-1}U = \bigcup_{k>0}\lbrace y \in \Pi^{-1}U : \alpha(y) > \frac{1}{k}, \rho(y) < k\rbrace,
$$
where the functions come from Theorem \ref{thm:NeilUnstableManif}, for $\eta < \chi_\mu$ a small positive constant.
There is therefore a $k_0 > 2$ such that the set
$$
B := \lbrace y \in \Pi^{-1}U : \alpha(y) > \frac{1}{k_0}, \rho(y) < k_0\rbrace
$$
is of positive $\overline{\mu}$ measure, i.e., $\overline{\mu}(B) > 0$, and of course $\Pi(B) \subset U$. Then there is an open interval $V \subset U$ such that $|V| < k_0^{-1}$ and $\omu(B \cap \Pi^{-1}V) > 0$. Set $V_B := B \cap \Pi^{-1}V$. By ergodicity, almost every point  $y \in Y$ enters $V_B$ infinitely often, at times $n_j(y)$, say. Then $\Pi y$ is contained in intervals $W_j$ of size less than 
$k_0 e^{-n_j(\chi-\eta)} $ mapped by $f^{n_j}$ with uniformly bounded distortion onto $V$. 
Taking $y \in \Pi^{-1} V$ and $n_j$ large, by the Intermediate Value Theorem there is a periodic point $p$ in $W_j \subset V \subset U$ and, by bounded distortion, it is repelling.

Similarly, almost every point in $I$ is contained in arbitrarily small intervals which get mapped onto $V$, so we conclude that the closure of the set of inverse images of $p$ has full measure.
%Let $y$ be such a point contained in $V_B$. Let $n$ be sufficiently large that $k_0
 %e^{-n_j(\chi-\eta)}< \dist(\Pi y, \partial V)$ and be such that $F^n y \in V_B$. Then 
%and let $V_y$ be given by Theorem \ref{thm:NeilUnstableManif}. Since $\alpha(y) > k_0^{-1}$, $\PiY_y \supset V$.  

%In particular, $\mu$ almost every point $x$ in $I$ is the projection of a point $y$ which enters $V_B$ at infinitely often at times $n_j$ say. But then, by Theorem \ref{thm:NeilUnstableManif}, $x$ is contained in an interval of size
%$k_0 e^{-n_j(\chi-\eta)}$ which is mapped by $f^{n_j}$ onto $V$ with uniformly bounded distortion. Thus almost every point is contained in arbitrarily small intervals mapped by iterates of $f$ onto $V$. Taking $x$ in the interior of $V$, 
\eprf

Now we shall prove a general lemma which we were unable to find in the literature. 

We say that a set $V$ is a one-sided neighbourhood of a point $p$ if it contains a small interval $(p, p+ \varepsilon)$ for some non-zero real $\varepsilon$.
\begin{lem} \label{lem:pwmono}
    Let $f: \bigcup_{j=1}^n U_j \to U$ be a piecewise continuous map, continuous on each (pairwise disjoint) open interval $U_j$, for which $f_{|U_J}$ extends to a continuous map on the compact interval $\overline{U_j}$ for each $j$. 

    Suppose $\mu$ is an ergodic, invariant, probability measure, and $W \subset f(W)$ is an interval of positive measure. Then there exists $k$ such that $\mu(f^k(W)) = 1$. 
\end{lem}
\beginpf
Let $Y$ denote the  boundary set $\partial (\bigcup_j U_j)$ which is finite. 
Let $W_j = \bigcup_{i=0}^j f^i(W) = f^j(W)$. Since $Y$ is finite, there is an $N$ such that for all $j > N$, $W_j \setminus W_N$ does not form a one-sided neighbourhood of any point in $Y$. 

Suppose $j \geq N$ and $V$ is a connected component of $W_j \setminus Y$. Suppose there is some  $k \geq 1$  such that $f^k(V) \cap W_N \ne \emptyset$ and let $l$ be the least such $k$. 
Then $f^l_{|V}$ is continuous. 

For each connected component $V$ of $W_N \setminus Y$, let $n_V$ be the minimal $k\geq 1$ such that $f^k(V) \cap W_N \ne \emptyset$ if such a $k$ exists; otherwise set $n_V = 0$. 
Let $M$ be the maximum of the $n_V$, noting that $W_N$ has a finite number of connected components. 
Let $X$ be a connected component of $W_j$ for some $j\geq M$. Suppose $X$ does not contain a connected component of $W_M$. Let $y \in X$ and let $l \geq 1$ be minimal such that $y \in f^l(V)$, where $V$ is some component of $W_N\setminus Y$. Then $l>M$. In particular, $n_V = 0$, so $f^i(y) \notin W_N$ for all $i\geq 0$. This holds for each $y \in X$. Therefore $f^i(X) \cap W_N = \emptyset$ for all $i \geq 0$. By ergodicity, $X$ then has null measure. 

Let $V_1, \ldots, V_r$ denote the connected components of $W_M$.
For $j\geq M$ and $1 \leq k \leq r$, let $V^j_k$ denote the connected component of $W_j$ containing $V_k$. These are the only connected components of $W_j$ which may have positive measure.

Let $W_\infty = \bigcup_{j\geq 0}  W_j$. 
We shall say that a point $x$ \emph{has the one-sided property} if 
 $W_\infty$ contains arbitrarily small, nested,  one-sided neighbourhoods $J_k$ of the $x$, that each $J_k$ has positive measure, and that $J_k \not\subset W_j$ for any $j < \infty$ and that $|J_k| \to 0$. 
 
 Note that if the conclusion of the lemma does not hold then there exists a point $x$ which has the one-sided property. Suppose this is so. We must arrive at a contradiction. 

 Points with the one-sided property belong to the finite set $\bigcup_{k=1}^r \partial \bigcup_{j\geq M} V^j_k$. 
If $y$ is in the interior of $W_\infty$, then $f^k(y)$ does not have the one-sided property for any $k\geq 0$. 
Suppose $x$ has the one-sided property and let $J_x$ be a corresponding one-sided neighbourhood of $x$. Then 
there is a sequence $\{y_n\}$  for which $f(y_n) \in J_x$ and converges to $x$ and for which each neighbourhood of each $y_n$ has positive measure. Replacing by a subsequence if necessary, we can assume the sequence $\{y_n\}$ is monotone. Let $y$ be its limit. 
Then $f( (y_n, y) )$ contains $(f(y_n), x) \subset J_x$. Since $x$ has the one-sided property, $(y_n,y) \not\subset W_j$ for any $j>0$, so $y$ has the one-sided property. Moreover, if $y_n \in U_i$ for all large $n$, then the continuous extension of $f$ to $\overline{U_i}$ maps $y$ to $x$. 

%Iterating this argument, since the number of points having the one-sided property is finite, it follows that there exists 
%an integer $k\geq 1$, a point $z$ having the one-sided property, a corresponding one-sided neighbourhood $J_z$, and a monotone sequence $\{z_n\}$ of points converging to $z$ in a one-sided neighbourhood $J_z$ of $z$ such that $f^k(z_n)$ converges to $z$ in $J_z$ and such that each neighbourhood of each $z_n$ has positive measure. 

It follows that (the finite collection of) points having the one-sided property are periodic, in the sense that if $x,J_x$ are as before and $J_x$ is sufficiently small then some iterate $f^k(J_x)$ is a small one-sided neighbourhood of $x$. 
We used here that on each $U_i$ the map $f$ extends continuously to the boundary.

The point $x$ is a fixed point of the continuous extension $g$ of $f^k_{|J_x}$ to the closure $\overline{J_x}$.
If $ g(J_x) \subset J_x$ for a sequence of arbitrarily small  $J_x$ then $x$ must be a periodic orbit supporting the measure, by ergodicity.
Suppose the measure is not supported on a periodic orbit (the lemma is trivial otherwise).  
Then $x$ is (one-sided, topologically-) repelling fixed point for $g$.
For mass to accumulate near $x$ in $J_x$, there must be a point $y \ne x$ and arbitrarily small, positive measure,  one-sided neighbourhoods $J_y$ of $y$ mapped by $f^k$ onto one-sided neighbourhoods of $x$ in $J_x$. 
But $y, J_y$ are not periodic in the above sense, so $y$ does not have the one-sided property. In particular, some sufficiently small $J_y$ is contained in some $W_j$. Then $W_{j+k} \supset f^k(J_y)$ contains some $J_x$, contradicting the one-sided property.  
\eprf

So far, our cusp maps may have an infinite (but countable) number of discontinuities. In the following definition and lemma, we wish to consider maps for which the number of serious discontinuities is finite. 
\begin{dfn} We say a cusp map $f: \bigcup_j I_j \to I$ has a \emph{finite number of discontinuities} if it extends to a continuous map $f_*$ on $\overline{\bigcup_j I_j} \setminus X$, where $X \subset I\setminus \bigcup_j I_j$ is a finite set of points called \emph{points of discontinuity}, and $X$ contains $\partial \overline{\bigcup_j I_j}$.
\end{dfn}

\begin{lem}\label{lem:exactness}
Let $f$ be a piecewise monotone \epsgamma  map with a finite number of discontinuities. Suppose  $\mu \in \scrM(f)$ has a positive finite Lyapunov exponent.

Then  given any interval $V$ of positive $\mu$-measure, there is a $j > 0$ such that $\mu\left(\bigcup_{k = 0}^j f^k(V)\right) =1$.
\end{lem}
\beginpf
This follows immediately from Lemmas \ref{lem:pwmono} and \ref{thm:periodicpointsdense}.
\eprf

\begin{dfn}[\cite{Graczyk:MetricAttractors}] An open interval $U$ is \emph{regularly returning} if $f^n(\partial U) \cap U = \emptyset$ for all $n > 0$. This is also called a \emph{nice interval} in the literature.
\end{dfn}
For a cusp map $f$,  
if $A$ is a (necessarily open) connected component of $f^{-n}(U)$ and $B$ is a connected component of $f^{-m}(U)$ with $m \geq n$, it is easy to check that either $A\cap B = \emptyset$ or $B \subset A$, so  inverse images of regularly returning intervals are either nested or disjoint. Indeed, suppose $x \in \partial A \cap B$. Then $f^n(x) \in \partial U$ (since $f$ may be discontinuous, one uses that $x \in B$ to know that $f^n$ is defined on a neighbourhood of $x$), but $f^m(x) \in U$, contradiction. 

\begin{lem} \label{lem:idvlonesided}
Let $m$ be a non-atomic measure on $I$. Let $Z$ be the set of all $x \in I$ such that there exists an open interval $L$ with $m(L) = 0$ and $x \in \partial L$. Then $m(Z) = 0$.
\end{lem}
\beginpf
Each open interval $L$ such that $m(L) = 0$ is contained in a maximal such interval $V$. All such maximal intervals are pairwise disjoint. Thus there are at most a countable number of such maximal intervals $V$. But $m$ is non-atomic so $m(\partial V)= 0$ and $Z$ is the countable union of sets of measure zero and thus of measure zero itself.
\eprf

\begin{prop} \label{prop:regretall}
Let $f$ be a \epsgamma map and suppose $\mu \in \scrM(f)$ is non-atomic and has positive finite Lyapunov exponent.

Then  $\mu$-almost every point is contained in arbitrarily small regularly returning open intervals, the boundaries of which are repelling periodic points.
\end{prop}
\beginpf
By Lemmas \ref{lem:idvlonesided} and  \ref{thm:periodicpointsdense}, almost every point $x$ is accumulated on both sides by repelling periodic points. Take  one arbitrarily close periodic point, not in the orbit of $x$, on each side of $x$ and consider the partition defined by the orbits of these two points. The interior of each  partition element, in particular the partition element which contains $x$, is regularly returning. The result follows, since the measure is non-atomic.
\eprf

\end{section}
\section{Generating partition and the \\Dynamical Volume Lemma} \label{sec:dvl}

Given a map $f$ and a partition $\scrP$ we denote by $\scrP_k$ the partition $\bigvee_{i=0}^k f^{-i}\scrP$, and by $\scrP_k(x)$ the partition element containing the point $x$.
Let $Y, \alpha, \rho$ be as per Theorem \ref{thm:NeilUnstableManif} and let $\gamma_1$ be as per Theorem~\ref{thm:Unstable2}.  
Recall that  finite partitions  into intervals of monotonicity are defined in Definition \ref{dfn:mono}.
\begin{prop}\label{prop:partition}
Let $f$ be a \epsgamma map and suppose $\mu \in \scrM(f)$ is non-atomic and has positive finite Lyapunov exponent $\chi_\mu$. Suppose there exists a finite partition $\scrQ$ into intervals of monotonicity.

There exist a regularly returning interval $U$,  constants $K, \varepsilon > 0$, a finite partition $\scrP$ and a set $X$ of full measure with the following properties:
\begin{itemize}
\item
$\omu(A) > 0$, where
$$
A := \{y \in Y : \Pi y \in U, \alpha(y) \geq |U|, 
%\mbox{ \,} 
\rho(y) < K, \gamma_1(y) < K \mbox{ and } \dist(\Pi y, \partial U) > 2\varepsilon |U|\};
$$
\item
$\scrP = \{U, I\setminus U\} \vee \scrQ$;
\item $\scrP$ is generating;
\item for each $x \in X$ there exists a strictly monotone increasing sequence $\{n_j\}$ such that
$$
f^{n_j} : \scrP_{n_j}(x) \to U
$$
is a diffeomorphism with distortion bounded by $\log 2$;
\item
$\dist(x,\partial \scrP_{n_j}(x)) > \varepsilon |\scrP_{n_j}(x)|$;
\item
$$\lim_{j\to \infty} \frac{j}{n_j} > 0 \mbox{ and } \lim_{j\to\infty} \frac{n_j}{n_{j+1}} = 1.
$$
\end{itemize}
\end{prop}
\beginpf
There exists a $K > 1$ such that
$$
B:= \{y \in Y: \alpha(y) > K^{-1} \mbox{ and } \rho(y) < K\}$$
has positive measure. By Proposition \ref{prop:regretall}, one can cover a set of full measure by a countable collection of regularly returning intervals of diameter less than $K^{-1}$. Let $U$ be one such interval such that
$$
A_0:= B \cap \Pi^{-1}U$$
 has positive measure. Now $\mu(\partial U) = 0$, so there exists a set $A \subset A_0$ verifying the first claimed property of the lemma for some $\varepsilon >0$.

By the Birkhoff Ergodic Theorem, almost every $y \in Y$ returns to $A$ at times $n_j(y)$ with $\lim_{j\to \infty} n_j/n_{j+1} = 1$ and asymptotic frequency $\lim_{j\to\infty} (j/n_j) = \omu (A)$. 
Almost every point $x \in I$ is the projection of such a point $y$. 
Define $\scrP$ as above. 
By
Theorem~\ref{thm:NeilUnstableManif},
 $f^{n_j}$ maps $\scrP_{n_j}(x)$ diffeomorphically onto $U$; the partition elements shrink to the point $x$ as $j\to \infty$, so $\scrP$ is generating. 
\eprf

%\remark This also gives the existence of induced Markov maps.

%\section{Dynamical Volume Lemma}\label{sec:dvl}
With this partition we can  now give a very short proof of the following Dynamical Volume Lemma. For maps with ``bounded  $p$-variation", this was proven in \cite{HofbauerRaith:HDmeasure}. Maps with unbounded derivative do not have bounded $p$-variation. We denote by $\HD(\mu)$ the Hausdorff dimension of a measure $\mu$.
\begin{prop}\label{prop:dvl}
Let $f$ be a \epsgamma map and suppose $\mu \in \scrM(f)$ has positive finite Lyapunov exponent $\chi_\mu$. Suppose there exists a finite partition $\scrQ$ into intervals of monotonicity.

Then for $\mu$-almost every $x$,
$$
\lim_{r\to0^+} \frac{\log \mu(B(x,r))}{\log r} = \frac{h_\mu}{\chi_\mu};
$$
in particular, $\HD(\mu) = h_\mu/\chi_\mu$.
\end{prop}
\beginpf The latter equality follows, by Frostman's Lemma, from the former, which we now prove. Note that if $\mu$ is atomic the proposition is trivial.
We write $\chi$ for $\chi_\mu$.
Let $(Y, F, \overline{\mu})$ be the natural extension. Let $\eta < \chi$ be a small positive constant. Let $\scrP = \{P_0 = U , P_1,\ldots, P_d\}$ be a finite generating partition of $I$, $\varepsilon$ the constant and $X$ the set of full measure given by Proposition~\ref{prop:partition}.

By ergodicity and the Shannon-McMillan-Breiman Theorem (\cite{Parry:Book} p.39), there is  a set $X' \subset X$ of full measure such that the pointwise Lyapunov exponent $\lim_{n\to \infty} (1/n) \log |Df^n(x)| = \chi$ and
$$
h_\mu = \lim_{n \to \infty} \frac{-1}{n}\log \mu(\scrP_n(x))$$
 for all $x \in X'$. Now fix $x \in X'$ and let $n_j$ be given by Proposition~\ref{prop:partition}.
%Thus
%$$
%h_\mu = \lim_{j \to \infty} \frac{-1}{n_j}\log \mu(U_j).$$

Set $r_j := \dist(x, \partial \scrP_{n_j}(x))$ and $R_j = |\scrP_{n_j}(x)|$. We have (again by Proposition~\ref{prop:partition}) that $r_j \geq \varepsilon R_j$. Continuing,
$$
\lim_{j \to \infty} \frac{1}{n_j}\log r_j =
\lim_{j \to \infty} \frac{1}{n_j}\log R_j =
-\chi$$
and, since $\lim_{j \to \infty} \frac{n_j}{n_{j+1}} = 1$,
$$
\lim_{j \to \infty} \frac{\log r_j}{\log r_{j+1}} =
\lim_{j \to \infty} \frac{\log R_j}{\log R_{j+1}} =
1.$$

Thus, if $r_j \geq r \geq r_{j+1}$,
$$
\frac{\log\mu(B(x,r))}{\log r} \geq \frac{\log\mu(\scrP_{n_j}(x))}{\log r_{j+1}}
= \frac{-1}{n_j}\log{\mu(\scrP_{n_j}(x))} \frac{-n_j }{\log r_{j+1}}
$$
and the right-hand side tends to $(h_\mu/\chi)$ as $j\to \infty$.

If $R_j \geq r \geq R_{j+1}$,
$$
\frac{\log\mu(B(x,r))}{\log r} \leq \frac{\log\mu(\scrP_{n_{j+1}}(x))}{\log R_{j}}
= \frac{-1}{n_{j+1}}\log{\mu(\scrP_{n_{j+1}}(x))} \frac{-n_{j+1}}{\log R_j}
$$
and the right-hand side tends to $(h_\mu/\chi)$ as $j\to \infty$.

As $r \to 0$ one has $j \to \infty$ so
we conclude that
$$
\lim_{r \to 0}
\frac{\log\mu(B(x,r))}{\log r}
= \frac{h_\mu}{\chi}$$
as required.
\eprf

\section{Existence of Pesin partition} \label{sec:pesin}
An analogous result to the following was proven in \cite{Ledrappier:AbsCnsInterval} in a more restrictive  setting, but the same proof works. However it is unnecessarily complicated. We provide a short proof of
the existence of Pesin's partition, taking advantage of the properties of regularly returning intervals. 
\begin{prop}\label{thm:pesin}
Let $f$ be a piecewise-monotone \epsgamma  map. Suppose  $\mu \in \scrM(f)$ has positive finite Lyapunov exponent and denote the natural extension $(Y, F, \overline{\mu})$. 
Let $\xi$ be the measurable partition of $Y$ defined by 
$$
\xi = \bigvee_{i=0}^\infty F^i(\Pi^{-1}\scrP),
$$
where $\scrP$ is a partition given by Proposition \ref{prop:partition}. Then $\xi$ has the following properties:
\begin{enumerate}
\item the partition $\xi$ is increasing by $F$, $F^{-1}\xi > \xi$, and generates;
\item entropy of $\mu$ is given by $h(\mu) = H(F^{-1}\xi/\xi)$;
\item for $\overline{\mu}$ almost every point $y$, 
for all $k \geq 0$, $\Pi$ maps $F^{-k}(\xi(y))$ injectively into  an interval of monotonicity of $f$,
 where $\xi(y)$ denotes the element of $\xi$ containing $y$;
\item for $\overline{\mu}$ almost every point $y$, 
$
0 < \int_{\xi(y)} \Delta(y,y') dy' < \infty$, where the integration is with respect to the natural Lebesgue measure (i.e. the pullback of Lebesgue measure by $\Pi_{|\xi(y)}$) on each element of $\xi$ and 
$$
\Delta(y, y') = \lim_{n\to \infty }
 \frac{Df^n(\Pi F^{-n} y)}{Df^n(\Pi F^{-n} y')};
$$
\item let $U$ and $A$ be given by Proposition \ref{prop:partition}, so $\omu(A) > 0$: for all $y \in A$,  $\Pi \xi(y) = U$.
\end{enumerate}
\end{prop}
\beginpf
Since $\scrP$ is finite and generating, $h(f, \scrP) = h_\mu$. Let $\zeta$ be the partition of $Y$ given by $\zeta := \{\Pi^{-1} P : P \in \scrP\}$ and let $\xi := \bigvee_{i=0}^\infty F^i \zeta$. Then $\zeta$ is a finite, generating partition of $Y$, and 
$$
h_{\overline{\mu}} = h(F, \zeta) = H(\zeta | \bigvee_{i > 0} F^i \zeta) = H(F^{-1}\xi | \xi) = h(f, \scrP) = h_\mu.
$$
That $\Pi$ maps $F^{-k}\xi(y)$ injectively into an interval of monotonicity holds since $\scrP$ is a refinement of the partition into intervals of monotonicity. 

For $y \in A$, that $\Pi \xi(y) = U$ follows from the regularly returning property. For almost every $y \in Y$, there is an $n\geq0$ such that $F^ny \in A$. Then $\xi(y) \supset F^{-n}(\xi(F^ny))$, so $\Pi \xi(y)$ contains an open interval.

It remains to show that the integral is positive and finite, which follows as per \cite{Ledrappier:AbsCnsInterval}: For clarity, let us write $[F^{-k} \xi]$ instead of $F^{-k}\xi$. For all $k \geq 0,$ one has $\xi(y) = F^k([F^{-k}\xi](F^{-k}y))$, so the 
projection of the partition element $[F^{-k}\xi](F^{-k}y)$ contains an open interval for almost every $y$. For almost every $y$, $|Df(\Pi F^{-i}y)|$ is positive and finite for all $i$, and by ergodicity   
there exists a $k \geq 0$ such that $F^{-k}y \in A$. One has 
$$
\int_{\xi(y)} \Delta(y, y') dy' = \int_{F^k([F^{-k}\xi](F^{-k}y))} \Delta(F^{-k}y, F^{-k}y') \prod_{i=1}^k \frac{Df(\Pi F^{-i} y)}{Df(\Pi F^{-i} y')} dy'  
$$
$$
= \prod_{i=1}^k |Df(\Pi F^{-i}y)| \int_{[F^{-k}\xi](F^{-k}y)} \Delta(F^{-k}y, y') dy'.
$$
By the distortion bound of Proposition \ref{prop:partition}, the last integrand is bounded inside $(2^{-1}, 2)$ since $[F^{-k}\xi](F^{-k}y) \subset \xi(F^{-k}y)$ and $F^{-k}y \in A$. Thus the integrals are positive and finite, completing the proof. 
\eprf
\section{Absolutely continuous measures}\label{sec:ledrappartition}
The Rohlin decomposition $p(y,\cdot)$ for the measure $\overline{\mu}$ with respect to the partition $\xi$ is a conditional probability measure on each partition element of $Y$ such that, for any measurable set $B \subset Y$ one has
$$
\overline{\mu}(B) = \int_Y p(y,B) d\overline{\mu} = \int_Y p(y, B\cap\xi(y)) d\overline{\mu}.
$$
By Proposition \ref{thm:pesin}, if $n > 0$,
\begin{equation} \label{eqn:nentropy}
nh_\mu = H(F^{-n}\xi / \xi) = -\int \log p(y, [F^{-n}\xi](y)) d\overline{\mu}.
\end{equation}

\begin{prop}\label{thm:rohlin}
Let $f$ be a piecewise-monotone \epsgamma  map. Suppose $\mu \in \scrM(f)$ has  Lyapunov exponent $\chi_\mu$ and entropy $h_\mu$ satisfying $0 < \chi_\mu = h_\mu < \infty$. Let $(Y, F, \overline{\mu})$ be the natural extension. 

Then the Rohlin decomposition for the measure $\overline{\mu}$ with respect to the partition $\xi$ of Proposition \ref{thm:pesin} is given by $q(y,B)$, for $y$ in $Y$ and $B$ a measurable subset, where
\begin{equation}\label{eqn:rohlin}
q(y, B) := \frac{\int_{B \cap \xi(y)} \Delta(y,y') dy'}{\int_{\xi(y)} \Delta(y,y') dy'}.
\end{equation}
\end{prop}
\beginpf
The proof carries over from \cite{Ledrappier:AbsCnsInterval}, proposition 3.6, without modification. 
\eprf

\begin{cor} \label{cor:cuspacim}
Let $f$ be a piecewise-monotone \epsgamma  map. Suppose $\mu \in \scrM(f)$ has positive finite Lyapunov exponent $\chi_\mu$.

If $h_\mu = \chi_\mu$, or equivalently if $\HD(\mu) = 1$, then $\mu$ is absolutely continuous.
\end{cor}
\beginpf Let $B$ be a subset of zero Lebesgue measure. Then $q(y,B) = 0$ for all $y$, so $\overline{\mu}(\Pi^{-1}B) = 0 = \mu(B)$.
\eprf

%The following result is new even for $C^{1+\epsilon}$ maps. Keller \cite{Keller:density} has proven it under the assumption that the critical points of the smooth map $f$ are non-flat. 
\begin{prop}\label{thm:cuspacim}
Let $f$ be a piecewise-monotone \epsgamma  map. Suppose $\mu \in \scrM(f)$ has positive finite Lyapunov exponent.

If $\mu$ is absolutely continuous with respect to Lebesgue measure  then there exist $\nu > 0$ and an open interval such that the density of $\mu$ is bounded from below by some constant $\nu > 0$ Lebesgue almost everywhere on the interval.
\end{prop}
\beginpf
Let $\xi, A, U$ be given by Proposition \ref{thm:pesin}. 
 For all $y \in A$, $\Pi \xi(y) = U$ and  $2^{-1} \leq \Delta(y,y') \leq 2$ if $y' \in \xi(y)$.
Thus, for $y \in A$, the density of the Rohlin decomposition on the partition element $\xi(y)$ containing $y$ is bounded inside $[4^{-1}|U|^{-1},4 |U|^{-1}]$.  
Let $A' := 
 \bigcup_{y\in A} \xi(y)$. If we set $\omu_A := \omu_{|A'}$, then  $\Pi_* \omu_A$ has support  $U$ and the density of $\Pi_* \omu_A$ is bounded inside 
$$[4^{-1}|U|^{-1}\omu(A'),4 |U|^{-1}\omu(A')].$$  
Since $\omu_A$ is a restriction of $\omu$, the density of $\Pi_* \omu_A$ is less than that of $\mu$ almost everywhere. In particular, the density of $\mu$ on $U$ is bounded from below by $\nu:= 4^{-1}|U|^{-1}\omu(A') > 0$ as required. 
\eprf
\begin{lem}\label{lem:InvariantDensity}
Let $g$ be a piecewise $C^1$ map with an absolutely continuous invariant probability measure with density $\rho$. Then
$$
\rho(x) = \sum_{w\in g^{-1}x} \frac{1}{|Dg(w)|}\rho(w).
$$
In particular, $\rho(g(x)) \geq \frac{1}{|Dg(x)|}\rho(x)$.
\end{lem}
\beginpf
This is just the change of variables formula.
\eprf

\begin{thm}\label{thm:densitybddbelow}
Let $f : I \to I$ be a piecewise-monotone \epsgamma  map with only a finite number of discontinuities. 
Suppose $\mu \in \scrM(f)$ has positive finite Lyapunov exponent and that $\mu$ is absolutely continuous with respect to Lebesgue measure.

Then the support of $\mu$ is a finite union of intervals $X$ 
on which $\mu$ is equivalent to Lebesgue. 

Moreover, on every compact subset of $X$ disjoint from the forward orbit of poles of $f$ the density is bounded away from 0. In particular, if $f$ has no poles the density is bounded away from 0 on $X$.
%if, in addition, $f$ has no poles, then the density of $\mu$ on $X$ is bounded from below by a positive constant $\delta > 0$.
\end{thm}
\beginpf
Let $U$ be the interval  given by  Proposition \ref{thm:cuspacim} on which the density is bounded away from 0. By Lemma \ref{lem:exactness}, there exists a $j > 0$ such that the closure $X$ of $\bigcup_{k=0}^j f^kU$ contains the support of $\mu$, and thus equals the support of $\mu$.
The result then follows from  Lemma \ref{lem:InvariantDensity}.
\eprf

\bigskip

\noindent \emph{Proof of Theorem \ref{thm:response}:}
%\begin{thm}\label{thm:response}
%Let $f: I \to I$ be a $C^{1+\epsilon}$ map of the interval $I$. Suppose  $\mu \in \scrM(f)$ has positive entropy, or  equivalently positive Lyapunov exponent, and that $\mu$ is absolutely continuous with respect to Lebesgue measure. 
%
%Then the support of $\mu$ is a finite union of intervals on which
%$$
%\int \log |Df(x)| dx > -\infty,
%$$
%where integration is with respect to Lebesgue measure. 
%\end{thm}
%\beginpf
%We remark by Ruelle's Inequality, positive entropy implies positive Lyapunov exponent, and the Dynamical Volume Lemma gives the other direction of the equivalence. 
We can apply Theorem \ref{thm:densitybddbelow}, since $f$ extends to be a continuous cusp map on some larger interval. The lower bound on the density implies that anything integrable with respect to $\mu$ is integrable with respect to Lebesgue measure on the support of $\mu$.
\eprf

\section{Induced Markov maps} \label{sec:markov}

When trying to prove the existence of absolutely continuous invariant probability measures, a standard and fruitful technique is to show the existence of an expanding induced map whose domain has full measure in its range. 
One can spread the absolutely-continuous, invariant, probability measure for the induced map to get an absolutely-continuous, invariant measure for the original map. 
If the return time for such a \emph{Markov map} is integrable with respect to Lebesgue measure, 
then the resultant measure is finite and can be normalised to give a probability measure. 

A natural question is whether for \emph{all} acips the measure can be produced from such an induced map. Henk Bruin \cite{Bruin:Markov0} has shown that 
this is the case for unimodal maps with negative Schwarzian derivative and non-flat critical points (this has now been extended to multimodal maps with negative Schwarzian derivative and non-flat critical points \cite{BT:EquilibriumInterval}).
 %He also shows that it is true when ``natural'' Markov maps exist, but does not prove  that they do outside of the unimodal negative Schwarzian non-flat critical point case. 
We prove a  stronger result, dropping the condition on non-flatness of critical points, admitting poles and weakening the condition on the number of singularities. 

\begin{dfn} Suppose $I' \subset I$ and $f : I' \to I$, where $I$ is an interval. 
Let $\{U_i\}$ be a finite or countable collection of disjoint open subintervals of an open interval $U \subset I'$.
We call a map $\phi : \bigcup_i U_i \to U$ an  \emph{expanding induced 
%$C^{1+\epsilon}$ 
Markov map} if 
\begin{itemize}
\item
$\phi$ restricted to each $U_i$ is a diffeomorphism onto $U$;
\item   $|D\phi| \geq \lambda > 1$ for some constant $\lambda$;
\item
there exist constants $C, \epsilon>0$ such that for each $i$, for all $x,x' \in U_i$,
$$
|D\phi(x) - D\phi(x')| \leq C |\phi(x) - \phi(x')|^\epsilon;
$$
\item 
there exists $\{n_i\}$ such that $\phi_{|U_i} = f^{n_i}_{|U_i}$.
\end{itemize}
If moreover $U \setminus \bigcup_i U_i$ has zero Lebesgue measure, then we call $\phi$ \emph{full}. 

Let $n(x) := n_i$ if $x \in U_i$. If $\phi$ is full and
$$
\int_U n(x) dx = \sum_i n_i |U_i| < \infty
$$
then we say $\phi$ has \emph{integrable return time}.
\end{dfn}
It follows easily from Proposition \ref{prop:partition} that if $f$ is a cusp map and $\mu \in\scrM(f)$ has positive finite Lyapunov exponent, then 
there exists an expanding induced Markov map such that $\mu(\bigcap_{j\geq0} \phi^{-j}(U)) = \mu(U) >0$ (note that holding of the H\"older condition comes from Theorem~\ref{thm:Unstable2} via the $\gamma_1$ of Proposition~\ref{prop:partition}). We want to show more than this in the case that $\mu$ is absolutely continuous.
 
The Folklore Theorem (see \eg \cite{DeMeloVanStrien}) implies that if $\phi$ is a full expanding induced Markov map then $\phi$ has a unique absolutely continuous invariant probability measure, $\nu$ say, whose density is bounded away from zero and infinity on $U$.

If $\phi$ has integrable return time then, as is well known,
$$
\sum_i \sum_{j=0}^{n_i-1}f^j_*\nu_{|U_i}
$$
is a finite, ergodic, absolutely-continuous, invariant measure for $f$ which, when normalised, is an ergodic acip $\mu$ for $f$. We say $\mu$ \emph{is generated by} $\phi$.

\begin{prop} \label{thm:cuspacimMarkov}
Let $f$ be a piecewise-monotone \epsgamma  map. 
Suppose $\mu \in \scrM(f)$ has positive finite Lyapunov exponent and that $\mu$ is absolutely continuous with respect to Lebesgue measure.

Then $\mu$ is generated by a full expanding induced Markov map for $f$. 
%If $f$ is $C^r$ then $\mu$ is generated by a $C^r$ induced Markov map for $f$.
\end{prop}
\beginpf
Let $A', U$ be defined as per the proof of Proposition \ref{thm:cuspacim}. For $y \in Y$ let $r_1(y) :=  \inf\{n\geq 1: F^ny \in A'\}$.
Inductively define $r_{k+1}(y)  := r_k(y) + r(F^{r_k(y)}y)$ for $k \geq 1$
and set $n_k(x) := \min\{r_k(y) : \Pi(y) = x\}$. 
These are defined on sets of full measure, in particular for Lebesgue almost every $x \in U$. Recall we showed that the density of $\Pi_* \omu_{|A'}$ is bounded from below on $U$ by some $\nu >0$. 

By an easy generalisation of Kac's Lemma, $\int_{A'} r_k(y) d\omu = k$ for each $k \geq 1$. Then 
$$
k \geq \int_U n_k(x) d \Pi_* \omu_{|A'} \geq \nu^{-1}\int_U n_k(x) dx.
$$
In particular $\int_U n_k(x) dx < \infty$.

Recall that for $y \in A'$, $\Pi \xi(y) = U$, and for $y \in A$, $\rho(y) < K$, so by Theorem \ref{thm:NeilUnstableManif}, $|Df^n| > (2K)^{-1} e^{n(\chi -\eta)}$ on $\Pi F^{-n}\xi(y)$. Choose $N$ such that  $(2K)^{-1} e^{N(\chi -\eta)} > 2$.

Let $D$ denote the set of $x\in U$ such that  $n_N(x)$ is defined. Since $n_N(x)$ is defined almost everywhere, the Lebesgue measure of $U \setminus D$ is zero. For each $x \in D$ there exists $y \in \Pi^{-1}(x)$ such that $F^{n_N(x)}y  \in A'$. Set 
$$
U_x := \Pi F^{-r_N(y)} \xi( F^{r_N(y)}), 
$$
and note that since $\xi(y') = \xi(y)$ for all $y' \in \xi(y)$ and since $U$ is regularly returning, $U_x$ is a connected component of $D$. in particular, $n_N(x')$ is defined and constant on $U_x$. Let $\{U_i\}$ be the collection of connected components of $D$ and $m_i := n_N(x)$ for some $x \in U_i$.
Then define $\phi : \bigcup U_i \to U$ by $\phi_{|U_i} := f^{m_i}$. This map is a full expanding induced Markov map.

Thus there is an ergodic invariant absolutely continuous invariant probability measure $\mu'$ generated by $\phi$. The support of $\mu'$ coincides with that of $\mu$ and both have positive density so by ergodicity they are equal and $\mu$ is generated by $\phi$.
\eprf

\begin{thm}
Let $f$ be a piecewise-monotone $C^{1+\epsilon}$  map. 

Then $f$ has an ergodic, absolutely continuous, invariant, probability measure with positive finite Lyapunov exponent if and only if there exists a full expanding induced Markov map with integrable return time.
\end{thm}
\beginpf
One direction is given by Proposition \ref{thm:cuspacimMarkov}. On the other hand, if there exists an induced Markov map with integrable return time, the measure generated by it will be an ergodic, absolutely continuous, invariant probability measure. The entropy of the measure is positive because it is non-invertible almost everywhere on the range of the Markov map. Then Ruelle's Inequality implies that the Lyapunov exponent is positive, and it is finite because the derivative is bounded.
\eprf

\remark There are induced Markov maps with integrable return time for cusp maps such that the generated measure has non-integrable Lyapunov exponent, see the section after next.

\section{From Markov extension to natural extension} \label{sec:posentropy}

This section builds on the work of Hofbauer (\cite{Hofbauer:IntrinsicErgodicityI}) and Keller (\cite{Keller:Liftability}) and looks at the relation between their \emph{Markov extension} or \emph{Hofbauer tower} and the natural extension. 

Consider a piecewise-monotone map $f : \bigcup_{j=1}^d I_j \to I$ of the interval $I$ defined on a collection of open intervals $I_j$ and an invariant, ergodic, probability measure $\mu$.  As is well known, Hofbauer showed that if the entropy $h_\mu >0$, then $\mu$ lifts to an ergodic, conservative, $\widehat{f}$-invariant probability measure $\widehat{\mu}$ for the corresponding Markov extension $\widehat{f} : \widehat{\bigcup_{j=1}^d I_j} \to \widehat{I}$. 

Bruin and Todd use this in the following way: There is some interval $\widehat{J}$ at some level in the tower which has positive measure. By conservativity, the domain of the first return map $\widehat{\phi}$ to $\widehat{J}$ has full measure. The topological structure of the tower implies that, when chosen appropriately, $\widehat{\phi}$ is a homeomorphism from each connected component of its domain onto $\widehat{J}$. This then also holds for all iterates of $\widehat{\phi}$. By conservativity, $\bigcap_{n\geq0} \widehat{\phi}^{-n} (\widehat{J})$ has the same measure as $\widehat{J}$. See section 3, and in particular subsection 3.2, of \cite{BT:EquilibriumInterval} for details.  

We can look at the natural extension of the Markov extension. Since it projects onto the Markov extension which then projects onto the original system, it is none other than the natural extension $(Y, F, \overline{\mu})$ of the original system. Let $\widehat{\pi}$ denote the projection 
from $Y$ to $\widehat{I}$, and, as usual, $\Pi$ the projection from $Y$ to $I$. 

Let $A \subset Y$ denote $\widehat{\pi}^{-1}(\widehat{J})$. Then $A$ has positive measure. Let $\Phi: A \to A$ denote the first return map to $A$ and set  $A_* := \bigcap_{j\geq 0} \Phi^{-j} (A)$;
again by conservativity, $A_*$ has positive measure. The first return map $\Phi$ projects down to the map $\widehat{\phi}$ on $\widehat{J}$, that is, $\widehat{\pi} \circ \Phi = \widehat{\phi}\circ \widehat{\pi}$ since $\widehat{\phi}$ is a first return map [note that the projection down to the original system in general is not a first return map].  

Let $y \in A_*$ and $y_j = \Phi^{-j}y$. Then $\widehat{\phi}^j$ maps a neighbourhood $\widehat{W_j} \subset \widehat{J}$ of $\widehat{\pi}(y_j)$ homeomorphically onto $\widehat{J}$, so there is some $k_j$ and a corresponding neighbourhood $W_j$ of $\Pi y_j \in I$ mapped by $f^{k_j}$  homeomorphically onto $J$, 
the projection of $\widehat{J}$, and $f^{k_j} \circ \Pi = \Pi \circ \Phi^j$ on $\Pi^{-1}W_j$. 

    Then for each $x$ in $J$, there is a corresponding point $y^x \in \Pi^{-1}x$ such that for all $j\geq 0$, $\Pi\Phi^{-j}(y^x) \in W_j$. In particular, $\Pi\Phi^{-j}(\bigcup_{x\in J} y^x) = W_j$. It follows that for each $k \geq 0$, $\Pi F^{-k}(\bigcup_{x\in J} y^x)$ is an interval mapped homeomorphically by $f^k$ onto $J$. 

    We remark that almost every point $y$ lands in $A_*$ at some point. Thus the following holds. 
    \begin{prop} \label{prop:fibres}
         Suppose $f$ is a piecewise-monotone, piecewise-continuous map and $\mu$ is an ergodic, invariant, probability measure with positive entropy $\mu$. There is a measurable function $\theta$ on $Y$ such that for almost every $y$, there exists a set $V_y \subset Y$ with the following properties: 
    \begin{itemize}
    \item
        $y \in V_y$, $\theta(y) >0$ and $\Pi V_y = B(\Pi y, \theta(y))$;
        \item
        for each $m > 0$, $f^m : \Pi F^{-m} V_y \to \Pi V_y$ is a homeomorphism (in particular it is onto).
        \end{itemize}
    \end{prop}

    Now, the collection $\scrQ := \{I_1, \ldots, I_d\}$ defines a partition of the domain of $f$ and $f$ is a homeomorphism on each $I_j$ (onto its image). 
    %Let $\scrQ_n := \bigvee_{i=0}^{n-1} f^{-i}(\scrQ)$, so 
    If $f^n$ is a homeomorphism from some interval $W$ onto its image, then $W$ is entirely contained in one element of the partition $\scrQ$. 
    Denote by $\scrQ(x)$ the partition element containing the point $x$. 

    By Proposition \ref{prop:fibres}, for almost every $y$ there is a set $V_y \ni y$ for which $\Pi F^{-n}(V_y)$ is mapped homeomorphically onto the open interval $\Pi V_y$ by $f^n$. Thus $\Pi F^{-n}(V_y)$ is entirely contained in some element of $\scrQ$.

    We now provide the remaining arguments needed to show Theorem \ref{thm:NeilUnstableManif}. Lemma 13 almost holds: one must consider the intervals $I_j$ one at a time. Lemma 16 holds unchanged. Lemma 17 is modified slightly taking into account Lemma 13:
    \begin{lem} 
        Let $B_n := B(\Pi F^{-n}y, 2^{-1} e^{-n3\delta/\epsilon})$. For all $n \geq n(y)$, for all $x,x' \in B_n \cap \scrQ(\Pi F^{-n}(y)) $,
$$
    \left| \log|Df(x)| - \log|Df(x')|\right| < e^{-n\delta}.$$
    \end{lem}

Lemma 18 becomes:
    \begin{lem}
    For $n \geq n(y), f(B_{n+1} \cap \scrQ(\Pi F^{-(n+1)}y) \supset B(\Pi F^{-n} y, 2 e^{-n(\chi-\delta)}) \cap f(\scrQ(\Pi F^{-(n+1)}y))$. 
    \end{lem}

Lemma 19 holds under the additional assumption that $V_n \subset f(\scrQ(\Pi F^{-(n+1)}))$. This holds if $V$ is a subset of $\Pi V_y$, where $V_y$ comes from Proposition \ref{prop:fibres}. The remainder of the proof follows through. 

\begin{section}{Conjugacies with the Chebyshev map} \label{sec:cheb}
The Chebyshev quadratic map  $x \mapsto 4x(1-x)$ is a very special example as it is smoothly conjugate on the interior of $[0,1]$ to a piecewise linear map, namely the full tent map $T :  [0,1] \to [0,1]$:
$$
T(x) = 
\left\{ \begin{array}{ll}
2x & \mbox{if $0 \leq x \leq \frac{1}{2}$};\\
2 -2x & \mbox{if $\frac{1}{2} < x \leq 1$}.\end{array} \right. 
$$
We shall use it to construct some examples of cusp maps.
In Theorem \ref{thm:response} and Corollary \ref{cor:acimaone}, we showed that for $C^{1+\epsilon}$ maps there are no absolutely continuous invariant measures with positive Lyapunov exponent if the critical points are too flat. Here we shall find unimodal maps with acips with positive entropy and \emph{very} flat critical points. This is strange, but a contradiction is avoided because these maps have singularities at the boundary and the logarithm of the derivative is non-integrable, so the Lyapunov exponent of the acip does not exist. 

Topological Chebyshev maps will also provide examples of maps whose measures of maximal entropy have infinite Lyapunov exponent. 

\begin{prop} \label{prop:examples}
Let $TC$ denote the class of cusp maps or restrictions of cusp maps which are topologically conjugate by a conjugacy $h$ to the full tent map $T$. Let $TC^\infty$ denote the subclass of such maps such that $h$ is a $C^\infty$ diffeomorphism on the interior of the interval. For $f \in TC$ let $\mu_f$ denote the pullback by $h$ of the acip (actually Lebesgue measure) for $T$. For $f \in TC$, $\mu_f$ is the measure of maximal entropy. If $f \in TC^\infty$ then $\mu_f$ is also absolutely continuous.

Then there exist
\begin{enumerate}
\item $f \in TC^\infty$ with poles at the boundary such that the Lyapunov exponent of $\mu_f$ is 2;
\item $f \in TC^\infty$ such that the Lyapunov exponent of $\mu_f$ does not exist;
\item $f \in TC^\infty$ such that $f$ has a smooth parabolic fixed point at the boundary and the Lyapunov exponent of $\mu_f$ is 2;
\item $f \in TC$ such that $\mu_f$ has positive infinite Lyapunov exponent.
\end{enumerate}
\end{prop}
For each $\alpha > 0$ let $h_\alpha$ be a $C^\infty$ homeomorphism of $[0,1]$ with the following properties:
\begin{itemize}
\item $|Dh_\alpha| > 0$ except at $0$ and at $1$;
\item on a neighbourhood of $0$ one has $h_\alpha(x) = e^{-x^{-\alpha}}$;
\item the graph of $h$ has a point of central symmetry at $(\frac{1}{2}, \frac{1}{2})$, in other words $h(1-x) = 1 - h(x)$.
\end{itemize}

For each $\alpha > 0$ define $g_\alpha :  [0,1] \to [0,1]$ by 
$$
g_\alpha(x) := h_\alpha \circ T \circ h_\alpha^{-1}.
$$
Then $g_\alpha$ is conjugate to the tent map, symmetric and  has non-zero, non-infinite derivative everywhere except for at $0, \frac{1}{2},$ and $1$ (figure \ref{fig:chebmaps}). Denote $\frac{1}{2}$ by $c$. Of course, $g(0) = 0, g(c) = 1, g(1) = 0$. Let us calculate the derivative near 0.
Firstly note that, near zero, $h^{-1}(x) = (-\log x)^{-\frac{1}{\alpha}}$, so
$$
Dh^{-1}(x) = \frac{1}{\alpha} (-\log x)^{\frac{-1 - \alpha}{\alpha}} \frac{1}{x}.
$$
Clearly $|DT| = 2$. Next, $Dh(x) = \alpha x^{-1 -\alpha} e^{-x^{-\alpha}}$. Putting these together using the chain rule, 
\begin{eqnarray*}
Dg_\alpha(x) &=& Dh(2 h^{-1}(x)) 2 Dh^{-1}(x) \\
&=& \frac{\alpha}{\alpha} 2^{-1 - \alpha} (-\log x)^{\frac{1 + \alpha}{\alpha}} e^{-2^{-\alpha}(-\log x)^1} 2 (\log x)^{\frac{-1 - \alpha}{\alpha}} \frac{1}{x} \\  
&=& 2^{-\alpha} x^{2^{-\alpha}} \frac{1}{x}\\
&=& 2^{-\alpha} x^{2^{-\alpha} -1}.
\end{eqnarray*}
One can double-check using the conjugacy that near zero, $g_\alpha(x) = x^{2^{-\alpha}}$. 
The derivative near $1$ is the same but with negative sign.
We shall see soon that there is a transition in behaviour at $\alpha = 1$ which corresponds to a cusp at $0$ of root type $ = \frac{1}{2}$.

Now let us look at the derivative near the critical point. Here $Dh_\alpha \ne 0$ and we can assume it is almost constant $ \approx Dh_\alpha(c) =: l^{-1}$, say. If $x$ is close to $c$ then $h_\alpha^{-1}(c)$ is close to $\frac{1}{2}$ too, so it gets mapped by $T$ to a point close to $1$. We get 
\begin{eqnarray*}
|Dg_\alpha(x + c)| &=& \left|Dh_\alpha(Th_\alpha^{-1}(x+c))2Dh_\alpha^{-1}(x+c)\right|\\
&\approx& Dh_\alpha(2l |x|) 2 l \\
&=& \alpha (2l|x|)^{-1 -\alpha} e^{-(2l|x|)^{-\alpha}}.
\end{eqnarray*}
We have seen in Corollary \ref{cor:acimaone}, that at $\alpha =1$ there is some sort of a transition in other maps with this critical behaviour (polynomial times $e^{-|y|^{-\alpha}}$).

Lebesgue measure is both an acip and a measure of maximal entropy for $T$. We can pull this back by $h_\alpha$ to get the measure of maximal entropy for $g_\alpha$. This measure $\mu$ is absolutely continuous since $h_\alpha$ is smooth and is given by $d\mu(x) = Dh_\alpha^{-1}(x) dx$.
Let us try integrating $\log |Dg_\alpha|$ near $0$:
$$
\int \log |Dg_\alpha(x)| Dh^{-1}(x) \, dx 
$$
\begin{eqnarray*}
&\approx &\int (\log 2^{-\alpha} + (2^{-\alpha}-1)\log x) \left(\frac{1}{\alpha}(-\log x)^{\frac{-1-\alpha}{\alpha}} \frac{1}{x}\right) \, dx \\
&=& \int \log 2^{-\alpha} Dh^{-1}(x) - (2^{-\alpha} -1)\frac{1}{\alpha} (-\log x)^{\frac{-1}{\alpha}} \frac{1}{x} \, dx \\
&=& C_0 + \int C_1 \frac{1}{x(-\log x)^{\frac{1}{\alpha}}} \, dx,
\end{eqnarray*}
where $C_0$ and $C_1$ are positive constants. This last integral, well known in France as an \emph{int\'egrale de Bertrand}, is finite if and only if $\alpha < 1$. 

Near $c$, since $Dh$ is approximately a positive constant, $\log |Dg_\alpha|$ will be integrable with respect to $\mu$ if and only if it integrable with respect to Lebesgue. It is easy to see that it is integrable with respect to Lebesgue if and only if $\alpha < 1$.

Thus $g_\alpha$ has an absolutely continuous invariant probability measure with \textbf{finite} Lyapunov exponent if $\alpha < 1$. In this case, since $h_\alpha$ is smooth everywhere bar the boundary, every periodic point $p$ of period $k$ satisfies $|Df^k(p)| = 2^k$, so  the Lyapunov exponent is equal to 2.  

If $\alpha \geq 1$ then $g_\alpha$ has an acip whose  \textbf{Lyapunov exponent does not exist}. At the critical point $\int \log|Dg_\alpha| d\mu$ is $-\infty$; at 0 it is $+\infty$. If we write $\underline{\chi}(x)$ and $\overline{\chi}(x)$ for the lower and upper Lyapunov exponents of $g_\alpha$ at $x$, then, almost everywhere,
$$
\underline{\chi}(x) \leq 2 \leq \overline{\chi}(x).
$$

\begin{figure}[htb]
%\center{\includegraphics{../../Images/chebmaps}}
\center{\includegraphics[totalheight=0.18\textheight]{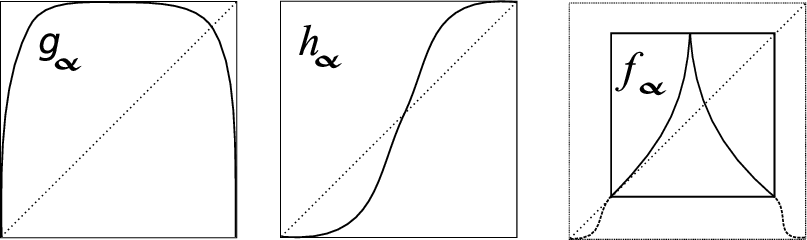}}
\caption{Graphs of $g_\alpha, h_\alpha$ and $f_\alpha$. \label{fig:chebmaps}}
\end{figure}

Now we conjugate in the other direction. 
For each $\alpha > 0$ define $f_\alpha :  [0,1] \to [0,1]$ by 
$$
f_\alpha(x) := h_\alpha^{-1} \circ T \circ h_\alpha.
$$
Again, $f_\alpha$ is symmetric, conjugate to the tent map and the derivative is non-zero and finite except possibly at $0, c, 1$. In fact the norm of the derivative is strictly positive and finite everywhere except for at $c$ (figure \ref{fig:chebmaps}). Let us calculate $Df_\alpha$ near $0$.
\begin{eqnarray*}
Df_\alpha(x) &=& Dh_\alpha^{-1}(2e^{-x^{-\alpha}}) 2 \alpha x^{-\alpha -1} e^{-x^{-\alpha}} \\
&=& \alpha^{-1} (-\log2 + x^{-\alpha})^{\frac{-1-\alpha}{\alpha}} \frac{1}{2} e^{x^{-\alpha}} 2 \alpha x^{-\alpha -1} e^{-x^{-\alpha}} \\
&=& (-\log 2 + x^{-\alpha})^{\frac{-1-\alpha}{\alpha}}x^{-\alpha - 1} \\
&=& (1 - x^\alpha \log 2)^{\frac{-1-\alpha}{\alpha}}.
\end{eqnarray*}
Thus $Df_\alpha(0) = 1$ and there is a parabolic fixed point at zero. It is straightforward to check that, for all $i$ such that $2 \leq i < \alpha + 1$, $D^if_\alpha(0)$ exists and $D^if_\alpha(0) = 0$. Also, $f_\alpha$ is $C^\infty$ near zero if and only if $\alpha>0$ is a natural number, in which case $D^{\alpha+1} f_\alpha(0) > 0$. 

Near the critical point we have
\begin{eqnarray*}
|Df_\alpha(c+x)| &=& \left| Dh_\alpha^{-1}(Th_\alpha(x+c))2Dh_\alpha(x+c)\right|\\
&\approx& Dh_\alpha^{-1}(2l^{-1} |x|) 2 l^{-1} \\
&=& 2l^{-1} \alpha^{-1} \left(-\log (l^{-1}2|x|)\right)^{\frac{-1-\alpha}{\alpha}} 2^{-1}l |x|^{-1}\\
&=& \alpha^{-1} |x|^{-1} \left(-\log (l^{-1}2|x|)\right)^{\frac{-1-\alpha}{\alpha}}. 
\end{eqnarray*}
Hence the derivative does indeed tend to infinity at $c$. Moreover, $\frac{1}{Df_\alpha}$ is $C^{\epsilon}$ near $c$, so an extension of $f_\alpha$ is an \epsgamma map. 

Pulling back Lebesgue measure, we get an absolutely continuous invariant probability measure $\mu$ for $f_\alpha$ with $C^\infty$ density $Dh(x)$, \textbf{despite the existence of a smooth parabolic point in the support of $\mu$}.
The Lyapunov exponent here is integrable and equal    to $2$.

Now we show that there is a cusp map $f$ in $TC$ whose measure of maximal entropy has infinite Lyapunov exponent.
Let $f : [0,1] \to [0,1]$ be a map in $TC$ with poles at the boundary and at the preimage $c$ of the turning point $\frac{1}{2}$, such that on $U$ a neighbourhood of zero, $f(x) = x^\alpha$ for some $\alpha$ with $0 < \alpha < 1$. See figure \ref{fig:topcheb}.

Let $f_1 := f_{|[0,c]}$ and let $p,N$ be such that $p = f_1^{-N}(c) \in U$ and set $p_n := f_1^{-n}(p)$. Then 
$$
p_n = p^{\frac{1}{\alpha^n}}.
$$
On $[p_{n+1}, p_n]$ one has $\alpha p_{n+1}^{\alpha -1} \leq Df \leq \alpha p_n^{\alpha - 1}$, and $\mu_f([p_{n+1}, p_n] = 2^{-(n+N+1)}$ since $f$ is conjugate to the full tent map $T_2$. So, on $[p_{n+1}, p_n]$, 
$$
\log \alpha + (\alpha-1)\frac{1}{\alpha^{n+1}}\log p \leq \log Df \leq 
\log \alpha + (\alpha-1)\frac{1}{\alpha^{n}}\log p.
$$
Subtracting $\log \alpha$ and integrating, 
\begin{eqnarray*}
(\alpha-1)\log p \sum_{i=0}^\infty \frac{1}{\alpha^{i+1}} 2^{-(i+N+1)} 
&\leq& - \mu_f([0,p])\log \alpha +  \int_0^p \log Df \, d\mu_f \\
&\leq& (\alpha-1)\log p \sum_{i=0}^\infty \frac{1}{\alpha^{i}} 2^{-(i+N+1)} 
\end{eqnarray*}
so, near zero, $\int \log Df \, d\mu_f$ is finite if $\alpha > \frac{1}{2}$ and positive infinite if $\alpha \leq \frac{1}{2}$. By definition, the derivative is bounded away from zero, so, if $\alpha \leq \frac{1}{2}$, then 
$$
\int_{[0,1]} \log |Df| \, d\mu = +\infty.
$$
\eprf
\begin{figure}[htb]
%\center{\includegraphics{../../Images/topcheb}}
\center{\includegraphics{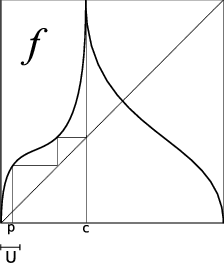}}
\caption{Graph of $f$. \label{fig:topcheb}}
\end{figure}

\end{section}

\section*{Acknowledgments}
The author would like to thank his former supervisor Jacek Graczyk for his encouragement and for many enlightening conversations. Most of the present work was carried out while a doctoral student at Université Paris-Sud. The author is currently a postdoctoral researcher at the University of Helsinki.
The author has benefited from funding from a Goldstine Fellowship (IBM), the Knut and Alice Wallenberg Foundation, the G\"oran Gustafsson Foundation and the EU Training Network CODY. Careful reading by a referee led to substantial improvement of this paper, for which the author is very grateful.

\bibliography{references}
\bibliographystyle{plain}

\end{document}